\documentclass[11pt,a4paper]{article}
\usepackage{mathrsfs}
\usepackage{indentfirst}
\usepackage{amsmath,amssymb,amsfonts,amsthm,graphics}
\input{epsf}
\usepackage{makeidx}
\usepackage{color}

\newtheorem{thm}{Theorem}
\newtheorem{cor}[thm]{Corollary}
\newtheorem{lem}{Lemma}

\newtheorem{obs}{Observation}

\usepackage{bbm}
\usepackage{mathrsfs}
\usepackage{amscd}
\usepackage{amsmath,amsfonts,amssymb,amscd}
\usepackage{indentfirst,graphics,epsfig,psfrag}
\input{epsf}
\usepackage{ifpdf}
\usepackage{enumerate}
\usepackage{appendix}
\usepackage{enumerate}
\usepackage{color}
\usepackage{lineno}

\makeatletter \@addtoreset{figure}{section} \makeatother
\makeatletter
\long\def\@makecaption#1#2{%
   \vskip 10\p@
   \setbox\@tempboxa\hbox{{#1}\ \ #2}%
   \ifdim \wd\@tempboxa >\hsize

       {#1}\ \ #2\par
   \else
       \hbox to\hsize{\hfil\box\@tempboxa\hfil}%
   \fi}
\makeatother

\def\qed{\hfill \rule{4pt}{7pt}}

\begin{document}
\title{\textbf{Graphs with large generalized
(edge-)connectivity}\footnote{Supported by NSFC No.11371205 and 11531011, and
PCSIRT.}}
\author{
\small  Xueliang Li$^1$, \ \ Yaping Mao$^{1,2}$\\[0.2cm]
\small $^1$Center for Combinatorics and LPMC-TJKLC\\
\small Nankai University, Tianjin 300071, China.\\[0.2cm]
\small $^2$Department of Mathematics, Qinghai Normal\\
\small University, Xining, Qinghai 810008, China\\[0.2cm]
\small Emails: lxl@nankai.edu.cn; maoyaping@ymail.com}
\date{}
\maketitle
\begin{abstract}
The generalized $k$-connectivity $\kappa_k(G)$ of a graph $G$,
introduced by Hager in 1985, is a nice generalization of the
classical connectivity. Recently, as a natural counterpart, we
proposed the concept of generalized $k$-edge-connectivity
$\lambda_k(G)$. In this paper, graphs of order $n$ such that
$\kappa_k(G)=n-\frac{k}{2}-1$ and $\lambda_k(G)=n-\frac{k}{2}-1$ for
even $k$ are characterized.\\[2mm]
\end{abstract}

\section{Introduction}

All graphs considered in this paper are undirected, finite and
simple. We refer to the book \cite{bondy} for graph theoretical
notation and terminology not described here. For a graph $G$, let
$V(G)$, $E(G)$, $\overline{G}$ denote the set of vertices, the set
of edges of $G$ and the complement, respectively. Let $d_G(v)$
denote the degree of the vertex $v$ in $G$. As usual, the
\emph{union} of two graphs $G$ and $H$ is the graph, denoted by
$G\cup H$, with vertex set $V(G)\cup V(H)$ and edge set $E(G)\cup
E(H)$. Let $mH$ be the disjoint union of $m$ copies of a graph $H$.
If $M$ is a subset of edges of a graph $G$, the subgraph of $G$
induced by $M$ is denoted by $G[M]$, and $G- M$ denotes the subgraph
obtained by deleting the edges of $M$ from $G$. If $M=\{e\}$, we
simply write $G-e$ for $G-\{e\}$. If $S\subseteq V(G)$, the subgraph
of $G$ induced by $S$ is denoted by $G[S]$. For $S\subseteq V(G)$,
we denote $G-S$ the subgraph obtained by deleting the vertices of
$S$ together with the edges incident with them from $G$. We denote
by $E_G[X,Y]$ the set of edges of $G$ with one end in $X$ and the
other end in $Y$. If $X=\{x\}$, we simply write $E_G[x,Y]$ for
$E_G[\{x\},Y]$. A subset $M$ of $E(G)$ is called a \emph{matching} of $G$
if the edges of $M$ satisfy that no two of them are adjacent in $G$.
A matching $M$ saturates a vertex $v$, or $v$ is said to be
\emph{$M$-saturated}, if some edge of $M$ is incident with $v$;
otherwise, $v$ is \emph{$M$-unsaturated}. If every vertex of $G$ is
$M$-saturated, the matching $M$ is \emph{perfect}. $M$ is a
\emph{maximum matching} if $G$ has no matching $M'$ with $|M'|>|M|$.

Connectivity and edge-connectivity are two of the most basic
concepts of graph-theoretic subjects, both in a combinatorial sense and an 
algorithmic sense. As we know, the classical connectivity has two
equivalent definitions. The \emph{connectivity} of a graph $G$,
written $\kappa(G)$, is the minimum size of a set $S\subseteq V(G)$
such that $G-S$ is disconnected or has only one vertex. If $G-S$ is
disconnected we call such a set $S$ a \emph{vertex cut-set} for $G$.
We call this definition the `cut' version definition of
connectivity. A well-known Menger's theorem provides an
equivalent definition of connectivity, which can be called the
`path' version definition of connectivity. For any two distinct
vertices $x$ and $y$ in $G$, the \emph{local connectivity}
$\kappa_{G}(x,y)$ is the maximum number of internally disjoint paths
connecting $x$ and $y$. Then
$\kappa(G)=\min\{\kappa_{G}(x,y)\,|\,x,y\in V(G),x\neq y\}$ is
defined to be the \emph{connectivity} of $G$. Similarly, the
classical edge-connectivity also has two equivalent definitions. The
\emph{edge-connectivity} of $G$, written $\lambda(G)$, is the
minimum size of an edge set $M\subseteq E(G)$ such that $G-M$ is
disconnected or has only one vertex. We call this definition the
`cut' version definition of edge-connectivity.
Menger's theorem also provides an equivalent definition of
edge-connectivity, which can be called the `path' version
definition. For any two distinct vertices $x$ and $y$ in $G$, the
\emph{local edge-connectivity} $\lambda_{G}(x,y)$ is the maximum
number of edge-disjoint paths connecting $x$ and $y$. Then
$\lambda(G)=\min\{\lambda_{G}(x,y)\,|\,x,y\in V(G),x\neq y\}$ is
defined to be the \emph{edge-connectivity} of $G$. For connectivity
and edge-connectivity, Oellermann gave a survey paper on this
subject, see \cite{Oellermann1}.

Although there are many elegant and powerful results on connectivity
in graph theory, the classical connectivity and edge-connectivity
also have their defects. So people want some generalizations of both
connectivity and edge-connectivity. For the `cut' version definition
of connectivity, we are looking for a minimum
vertex-cut with no consideration about the number of components of
$G-S$. Two graphs with the same connectivity may have differing
degrees of vulnerability in the sense that the deletion of a vertex
cut-set of minimum cardinality from one graph may produce a graph
with considerably more components than in the case of the other
graph. For example, the star $K_{1,n}$ and the path $P_{n+1}\ (n\geq
3)$ are both trees of order $n+1$ and therefore connectivity $1$,
but the deletion of a cut-vertex from $K_{1,n}$ produces a graph
with $n$ components while the deletion of a cut-vertex from
$P_{n+1}$ produces only two components. Chartrand et al.
\cite{Chartrand1} generalized the `cut' version definition of
connectivity. For an integer $k \ (k\geq 2)$ and a graph $G$ of
order $n \ (n\geq k)$, the \emph{$k$-connectivity} $\kappa'_k(G)$ is
the smallest number of vertices whose removal from $G$ produces a
graph with at least $k$ components or a graph with fewer than $k$
vertices. Thus, for $k=2$, $\kappa'_2(G)=\kappa(G)$. For more
details about $k$-connectivity, we refer to \cite{Chartrand1, Day,
Oellermann2, Oellermann3}. The $k$-edge-connectivity, which is a
generalization of the `cut' version definition of classical
edge-connectivity was initially introduced by Boesch and Chen
\cite{Boesch} and subsequently studied by Goldsmith in
\cite{Goldsmith1, Goldsmith2} and Goldsmith et al.
\cite{Goldsmith3}. For more details, we refer to \cite{Beineke,
Oellermann1}.

The generalized connectivity of a graph $G$, introduced by Hager
\cite{Hager}, is a natural and nice generalization of the `path'
version definition of connectivity. For a graph $G=(V,E)$ and a set
$S\subseteq V$ of at least two vertices, \emph{an $S$-Steiner tree}
or \emph{a Steiner tree connecting $S$} (or simply, \emph{an
$S$-tree}) is a subgraph $T=(V',E')$ of $G$ that is a tree with
$S\subseteq V'$. Two Steiner trees $T$ and $T'$ connecting $S$ are
said to be \emph{internally disjoint} if $E(T)\cap
E(T')=\varnothing$ and $V(T)\cap V(T')=S$. For $S\subseteq V(G)$ and
$|S|\geq 2$, the \emph{generalized local connectivity} $\kappa(S)$
is the maximum number of internally disjoint Steiner trees
connecting $S$ in $G$. Note that when $|S|=2$ a minimal Steiner tree
connecting $S$ is just a path connecting the two vertices of
$S$. For an integer $k$ with $2\leq k\leq n$, \emph{generalized
$k$-connectivity} (or \emph{$k$-tree-connectivity}) is defined as
$\kappa_k(G)=\min\{\kappa(S)\,|\,S\subseteq V(G),|S|=k\}$. Clearly,
when $|S|=2$, $\kappa_2(G)$ is nothing new but the connectivity
$\kappa(G)$ of $G$, that is, $\kappa_2(G)=\kappa(G)$, which is the
reason why one addresses $\kappa_k(G)$ as the generalized
connectivity of $G$. By convention, for a connected graph $G$ with
less than $k$ vertices, we set $\kappa_k(G)=1$. Set $\kappa_k(G)=0$
when $G$ is disconnected. This concept appears to have been
introduced by Hager in \cite{Hager}. It is also studied in
\cite{Chartrand2} for example, where the exact value of the
generalized $k$-connectivity of complete graphs are obtained. Note
that the generalized $k$-connectivity and the $k$-connectivity of a
graph are indeed different. Take for example, the graph $H_1$
obtained from a triangle with vertex set $\{v_1,v_2,v_3\}$ by adding
three new vertices $u_1,u_2,u_3$ and joining $v_i$ to $u_i$ by an
edge for $1 \leq i\leq 3$. Then $\kappa_3(H_1)=1$ but
$\kappa'_3(H_1)=2$. There are many results on the generalized
connectivity or tree-connectivity, we refer to \cite{Chartrand2,
LLMS, LLMY, LLSun, LLL2, LL, LLZ, LM, LM2, LM3, LM4, Okamoto}. Apart
from the concept of tree-connectivity, Hager also introduced another
tree-connectivity parameter, called the {\it pendant
tree-connectivity} of a graph in \cite{Hager}. For the
tree-connectivity, we only search for edge-disjoint trees which
include $S$ and are vertex-disjoint with the exception of the
vertices in $S$. But pendant tree-connectivity further requires the
degree of each vertex of $S$ in a Steiner tree connecting $S$ equal
to one. Note that it is a special case of the tree-connectivity.

As a natural counterpart of the generalized connectivity, we
introduced in \cite{LMS} the concept of generalized
edge-connectivity, which is a generalization of the `path' version
definition of edge-connectivity. For $S\subseteq V(G)$ and $|S|\geq
2$, the \emph{generalized local edge-connectivity} $\lambda(S)$ is
the maximum number of edge-disjoint Steiner trees connecting $S$ in
$G$. For an integer $k$ with $2\leq k\leq n$, the \emph{generalized
$k$-edge-connectivity} $\lambda_k(G)$ of $G$ is then defined as
$\lambda_k(G)= \min\{\lambda(S)\,|\,S\subseteq V(G) \ and \
|S|=k\}$. It is also clear that when $|S|=2$, $\lambda_2(G)$ is
nothing new but the standard edge-connectivity $\lambda(G)$ of $G$,
that is, $\lambda_2(G)=\lambda(G)$, which is the reason why we
address $\lambda_k(G)$ as the generalized edge-connectivity of $G$.
Also set $\lambda_k(G)=0$ when $G$ is disconnected. Results on the
generalized edge-connectivity can be found in \cite{LM, LM2, LMS}.

In fact, Mader \cite{Mader3} was studying an extension of Menger's
theorem to independent sets of three or more vertices. We know from
Menger's theorem that if $S=\{u,v\}$ is a set of two independent
vertices in a graph $G$, then the maximum number of internally
disjoint $u$-$v$ paths in $G$ equals the minimum number of vertices
that separate $u$ and $v$. For a set $S=\{u_1,u_2,\cdots,u_k\}$ of
$k$ vertices $(k\geq 2)$ in a graph $G$, an \emph{$S$-path} is
defined as a path between a pair of vertices of $S$ that contains no
other vertices of $S$. Two $S$-paths $P_1$ and $P_2$ are said to be
\emph{internally disjoint} if they are vertex-disjoint except for
their endvertices. If $S$ is a set of independent vertices of a
graph $G$, then a vertex set $U\subseteq V(G)$ with $U\cap
S=\varnothing$ is said to \emph{totally separate $S$} if every two
vertices of $S$ belong to different components of $G- U$. Let $S$ be
a set of at least three independent vertices in a graph $G$. Let
$\mu(G)$ denote the maximum number of internally disjoint $S$-paths
and $\mu'(G)$ the minimum number of vertices that totally separate
$S$. A natural extension of Menger' s theorem may well be suggested,
namely: If $S$ is a set of independent vertices of a graph $G$ and
$|S|\geq 3$, then $\mu(S)=\mu'(S)$. However, the statement is not
true in general. Take the above graph $H_1$ for example. For
$S=\{v_1,v_2,v_3\}$, $\mu(S)=1$ but $\mu'(S)=2$. Mader proved that
$\mu(S)\geq \frac{1}{2}\mu'(S)$. Moreover, the bound is sharp.
Lov\'{a}sz conjectured an edge analogue of this result and Mader
proved this conjecture and established its sharpness. For more
details, we refer to \cite{Mader3, Mader4, Oellermann1}.

In addition to being natural combinatorial measures, the Steiner
Tree Packing Problem and the generalized edge-connectivity can be
motivated by their interesting interpretation in practice as well as
theoretical consideration. From a theoretical perspective, both
extremes of this problem are fundamental theorems in combinatorics.
One extreme of the problem is when we have two terminals. In this
case internally (edge-)disjoint trees are just internally
(edge-)disjoint paths between the two terminals, and so the problem
becomes the well-known Menger theorem. The other extreme is when all
the vertices are terminals. In this case internally disjoint Steiner
trees and edge-disjoint trees are just edge-disjoint spanning trees
of the graph, and so the problem becomes the classical
Nash-Williams-Tutte theorem.

\begin{thm}{\upshape (Nash-Williams \cite{Nash}, Tutte \cite{Tutte})}\label{th1}
A multigraph $G$ contains a system of $\ell$ edge-disjoint spanning
trees if and only if
$$
\|G/\mathscr{P}\|\geq \ell(|\mathscr{P}|-1)
$$
holds for every partition $\mathscr{P}$ of $V(G)$, where
$\|G/\mathscr{P}\|$ denotes the number of crossing edges in $G$,
i.e., edges between distinct parts of $\mathscr{P}$.
\end{thm}

The generalized edge-connectivity is related to an important
problem, which is called the \emph{Steiner Tree Packing Problem}.
For a given graph $G$ and $S\subseteq V(G)$, this problem asks to
find a set of maximum number of edge-disjoint Steiner trees
connecting $S$ in $G$. One can see that the Steiner Tree Packing
Problem studies local properties of graphs, but the generalized
edge-connectivity focuses on global properties of graphs. The
generalized edge-connectivity and the Steiner Tree Packing Problem
have applications in $VLSI$ circuit design, see \cite{Grotschel1,
Grotschel2, Sherwani}. In this application, a Steiner tree is needed
to share an electronic signal by a set of terminal nodes. Another
application, which is our primary focus, arises in the Internet
Domain. Imagine that a given graph $G$ represents a network. We
choose arbitrary $k$ vertices as nodes. Suppose that one of the
nodes in $G$ is a \emph{broadcaster}, and all the other nodes are
either \emph{users} or \emph{routers} (also called \emph{switches}).
The broadcaster wants to broadcast as many streams of movies as
possible, so that the users have the maximum number of choices. Each
stream of movie is broadcasted via a tree connecting all the users
and the broadcaster. So, in essence we need to find the maximum
number of Steiner trees connecting all the users and the
broadcaster, namely, we want to get $\lambda (S)$, where $S$ is the
set of the $k$ nodes. Clearly, it is a Steiner tree packing problem.
Furthermore, if we want to know whether for any $k$ nodes the
network $G$ has the above properties, then we need to compute
$\lambda_k(G)=\min\{\lambda (S)\}$ in order to prescribe the
reliability and the security of the network.

The following two observations are easily seen from the definitions.
\begin{obs}\label{obs1}
Let $k,n$ be two integers with $3\leq k\leq n$. For a connected
graph $G$ of order $n$, $\kappa_k(G)\leq \lambda_k(G)\leq
\delta(G)$.
\end{obs}
\begin{obs}\label{obs2}
Let $k,n$ be two integers with $3\leq k\leq n$. If $H$ is a spanning
subgraph of $G$ of order $n$, then $\lambda_k(H)\leq \lambda_k(G)$.
\end{obs}

Chartrand et al. in \cite{Chartrand2} got the exact value of the
generalized $k$-connectivity for the complete graph $K_n$.

\begin{lem}{\upshape \cite{Chartrand2}}\label{lem1}
For every two integers $n$ and $k$ with $2\leq k\leq n$,
$\kappa_k(K_n)=n-\lceil k/2\rceil.$
\end{lem}

In \cite{LMS} we obtained some results on the generalized
$k$-edge-connectivity. The following results are restated, which
will be used later.

\begin{lem}{\upshape \cite{LMS}}\label{lem2}
For every two integers $n$ and $k$ with $2\leq k\leq n$,
$\lambda_k(K_n)=n-\lceil k/2\rceil.$
\end{lem}

\begin{lem}{\upshape \cite{LMS}}\label{lem3}
Let $k,n$ be two integers with $3\leq k\leq n$. For a connected
graph $G$ of order $n$, $1\leq \kappa_k(G)\leq \lambda_k(G)\leq
n-\lceil k/2 \rceil$. Moreover, the upper and lower bounds are
sharp.
\end{lem}

We also characterized graphs attaining the upper bound and obtained
the following result.

\begin{lem}{\upshape \cite{LMS}}\label{lem4}
Let $k,n$ be two integers with $3\leq k\leq n$. For a connected
graph $G$ of order $n$, $\kappa_k(G)=n-\lceil\frac{k}{2}\rceil$ or
$\lambda_k(G)=n-\lceil\frac{k}{2}\rceil$ if and only if $G=K_n$ for
even $k$; $G=K_n- M$ for odd $k$, where $M$ is a set of edges such
that $0\leq |M|\leq \frac{k-1}{2}$.
\end{lem}

One may notice that the graphs with
$\kappa_k(G)=n-\lceil\frac{k}{2}\rceil$ are the same as the graphs
with $\lambda_k(G)=n-\lceil\frac{k}{2}\rceil$. Our motivation of
this paper is to ask whether the graphs with
$\kappa_k(G)=n-\lceil\frac{k}{2}\rceil-1$ are different from the
graphs with $\lambda_k(G)=n-\lceil\frac{k}{2}\rceil-1$. In this
paper, graphs of order $n$ such that
$\kappa_k(G)=n-\lceil\frac{k}{2}\rceil-1$ and
$\lambda_k(G)=n-\lceil\frac{k}{2}\rceil-1$ for any even $k$ are
characterized.

\begin{thm}\label{th2}
Let $n$ and $k$ be two integers such that $k$ is even and $4\leq
k\leq n$, and $G$ be a connected graph of order $n$. Then
$\kappa_k(G)=n-\frac{k}{2}-1$ if and only if $G=K_n- M$ where $M$ is
a set of edges such that $1\leq \Delta(K_n[M])\leq \frac{k}{2}$ and
$1\leq |M|\leq k-1$.
\end{thm}

The above result can also be established for the generalized
$k$-edge-connectivity, which is stated as follows.

\begin{thm}\label{th3}
Let $n$ and $k$ be two integers such that $k$ is even and $4\leq
k\leq n$, and $G$ be a connected graph of order $n$. Then
$\lambda_k(G)=n-\frac{k}{2}-1$ if and only if $G=K_n- M$ where $M$
is a set of edges satisfying one of the following conditions:

$(1)$ $\Delta(K_n[M])=1$ and $1\leq |M|\leq
\lfloor\frac{n}{2}\rfloor$;

$(2)$ $2\leq \Delta(K_n[M])\leq \frac{k}{2}$ and $1\leq |M|\leq
k-1$.
\end{thm}

\section{Main result}

To begin with, we give the following lemmas.

\begin{lem}\label{lem5}
If $G$ is a graph obtained from the complete graph $K_n$ by deleting
a set of edges $M$ such that $\Delta(K_n[M])\geq r$, then
$\lambda_k(G)\leq n-1-r$.
\end{lem}
\begin{pf}
Since $\Delta(K_n[M])\geq r$, there exists at least one vertex, say
$v$, such that $d_{K_n[M]}(v)\geq r$. Then
$d_G(v)=n-1-d_{K_n[M]}(v)\leq n-1-r$. So $\delta(G)\leq d_G(v)\leq
n-1-r$. From Observation \ref{obs1}, $\lambda_k(G)\leq \delta(G)\leq
n-1-r$.
\end{pf}

\begin{cor}\label{cor1}
For every two integers $n$ and $k$ with $4\leq k\leq n$, if $k$ is
even and $M$ is a set of edges in the complete graph $K_n$ such that
$\Delta(K_n[M])\geq \frac{k}{2}+1$, then $\kappa_k(K_n- M)\leq
\lambda_k(K_n- M)< n-\frac{k}{2}-1$.
\end{cor}

\noindent \textbf{Remark 1.} From Corollary \ref{cor1}, if
$\kappa_k(K_n- M)=n-\frac{k}{2}-1$ or $\lambda_k(K_n-
M)=n-\frac{k}{2}-1$ for $k$ even, then $\Delta(K_n[M])\leq
\frac{k}{2}$.

In \cite{LMS}, we stated a useful lemma for general $k$.

Let $S\subseteq V(G)$ be such that $|S|=k$, and $\mathscr{T}$ be a
maximum set of edge-disjoint $S$-Steiner trees in $G$. Let
$\mathscr{T}_1$ be the set of trees in $\mathscr{T}$ whose edges
belong to $E(G[S])$, and $\mathscr{T}_2$ be the set of $S$-Steiner
trees containing at least one edge of $E_G[S,\bar{S}]$, where
$\bar{S}=V(G)- S$. Thus, $\mathscr{T}=\mathscr{T}_1\cup
\mathscr{T}_2$ (Throughout this paper, $\mathscr{T}$,
$\mathscr{T}_1$, $\mathscr{T}_2$ are defined in this way).

\begin{lem}{\upshape\cite{LMS}} \label{lem6}
Let $G$ be a connected graph of order $n$, and $S\subseteq V(G)$
with $|S|=k \ (3\leq k\leq n)$ and let $T$ be a $S$-Steiner tree. If
$T\in \mathscr{T}_1$, then $T$ contains exactly $k-1$ edges of
$E(G[S])$. If $T\in \mathscr{T}_2$, then $T$ contains at least $k$
edges of $E(G[S])\cup E_G[S,\bar{S}]$.
\end{lem}

\begin{lem}\label{lem7}
For every two integers $n$ and $k$ with $4\leq k\leq n$, if $k$ is
even and $M$ is a set of edges of the complete graph $K_n$ such that
$|M|\geq k$ and $\Delta(K_n[M])\geq 2$, then $\lambda_k(K_n- M)<
n-\frac{k}{2}-1$.
\end{lem}
\begin{pf}
Set $G=K_n- M$. We claim that there is an $S\subseteq V(G)$ with
$|S|=k$ such that $|M\cap \big(E(K_n[S])\cup
E_{K_n}[S,\bar{S}])|\geq k$ and $|M\cap \big(E(K_n[S])|\geq 1$.
Choose a subset $M'$ of $M$ such that $|M'|=k$. Suppose that
$K_n[M']$ contains $s$ independent edges and $r$ connected
components $C_1,\cdots,C_r$ such that $\Delta(C_i)\geq 2 \ (1\leq
i\leq r)$. Set $|V(C_i)|=n_i$ and $|E(C_i)|=m_i$. Then $m_i\geq
n_i-1$. For each $C_i \ (1\leq i\leq r)$, we select one of the
vertices having maximum degree, say $u_i$. Set $X_i=V(C_i)-u_i$.

If there exists some $X_j$ such that $|E(K_n[X_j])|\geq 1$, then we
choose $X_i\subseteq S$ for all $1\leq i\leq r$. Since
$|V(C_i)|=n_i$ and $X_i=V(C_i)- u_i$, we have $|X_i|=n_i-1$. By such
a choosing, the number of the vertices belonging to $S$ is
$\sum_{i=1}^r|X_i|=\sum_{i=1}^r(n_i-1)\leq \sum_{i=1}^rm_i\leq k-s$.
In addition, we select one endvertex of each independent edge into
$S$. Till now, the total number of the vertices belonging to $S$ is
$\sum_{i=1}^r|X_i|+s \leq (k-s)+s=k$. Note that if
$\sum_{i=1}^r|X_i|+s< k$, then we can add some other vertices in $G$
into $S$ such that $|S|=k$. Thus all edges of $E(C_i)$ and the $s$
independent edges are put into $E(K_n[S])\cup E_{K_n}[S,\bar{S}]$,
that is, all edges of $M'$ belong to $E(K_n[S])\cup
E_{K_n}[S,\bar{S}]$. So $|M\cap \big(E(K_n[S])\cup
E_{K_n}[S,\bar{S}])|\geq k$, as desired. Since $|E(K_n[X_j])|\geq
1$, it follows that $|M\cap \big(E(K_n[S])|\geq 1$, as desired.

Suppose that $|E(K_n[X_i])|=0$ for all $1\leq i\leq r$. Then each
$C_i$ must be a star such that $|E(C_i)|\geq 2$. Recall that $u_i$
is one of the vertices having maximum degree in $C_i$. Select one
vertex from $V(C_i)-u_i$, say $v_i$. Put all the vertices of
$Y_i=V(C_i)-v_i$ into $S$, that is, $Y_i\subseteq S$. Thus
$|Y_i|=n_i-1$. In addition, we choose one endvertex of each
independent edge into $S$. By such a choosing, the total number of
the vertices belonging to $S$ is
$\sum_{i=1}^r|Y_i|+s=\sum_{i=1}^r(n_i-1)+s\leq \sum_{i=1}^rm_i+s\leq
(k-s)+s=k$. Note that if $\sum_{i=1}^r|X_i|+s< k$ then we can add
some other vertices in $G$ into $S$ such that $|S|=k$. Thus all
edges of $E(C_i)$ and the $s$ independent edges are put into
$E(K_n[S])\cup E_{K_n}[S,\bar{S}]$, that is, and all edges of $M'$
belong to $E(K_n[S])\cup E_{K_n}[S,\bar{S}]$. So $|M\cap
\big(E(K_n[S])\cup E_{K_n}[S,\bar{S}])|\geq k$, as desired. Since
$|E(C_i)|\geq 2$, it follows that there is an edge $u_iw_i\in M\cap
K_n[S]$ where $w_i\in V(C_i)-\{u_i,v_i\}$, which implies that
$|M\cap \big(E(K_n[S])|\geq 1$, as desired.

From the above arguments, we conclude that there exists an
$S\subseteq V(G)$ with $|S|=k$ such that $|M\cap \big(E(K_n[S])\cup
E_{K_n}[S,\bar{S}])|\geq k$ and $|M\cap \big(E(K_n[S])|\geq 1$.
Since each tree $T\in \mathscr{T}_1$ uses $k-1$ edges in
$E(G[S])\cup E_G[S,\bar{S}]$, it follows that $|\mathscr{T}_1|\leq
({{k}\choose{2}}-1)/(k-1)=\frac{k}{2}-\frac{1}{k-1}$, which results
in $|\mathscr{T}_1|\leq \frac{k}{2}-1$ since $|\mathscr{T}_1|$ is an
integer. From Lemma \ref{lem6}, each tree $T\in \mathscr{T}_2$ uses
at least $k$ edges of $E(G[S])\cup E_G[S,\bar{S}]$. Thus
$|\mathscr{T}_1|(k-1)+|\mathscr{T}_2|k\leq
|E(G[S])|+|E_G[S,\bar{S}]|$, that is,
$|\mathscr{T}_1|k+|\mathscr{T}_2|k\leq
|\mathscr{T}_1|+{{k}\choose{2}}+k(n-k)-k$. So
$\lambda_k(G)=|\mathscr{T}|=|\mathscr{T}_1|+|\mathscr{T}_2|\leq
n-\frac{k}{2}-1-\frac{1}{k}<n-\frac{k}{2}-1$.\qed
\end{pf}\\

\noindent \textbf{Remark 2.} From Lemmas \ref{lem4} and \ref{lem7},
if $\kappa_k(K_n- M)=n-\frac{k}{2}-1$ or $\lambda_k(K_n-
M)=n-\frac{k}{2}-1$ for $k$ even and $2\leq \Delta(K_n[M])\leq
\frac{k}{2}$, then $1\leq |M|\leq k-1$, where $M\subseteq E(K_n)$.

\begin{lem}\label{lem8}
For every two integers $n$ and $k$ with $4\leq k\leq n$, if $k$ is
even and $M$ is a set of edges in the complete graph $K_n$ such that
$|M|\geq k$ and $\Delta(K_n[M])=1$, then $\kappa_k(K_n- M)<
n-\frac{k}{2}-1$.
\end{lem}
\begin{pf}
Let $G=K_n- M$. Since $\Delta(K_n[M])=1$, it follows that $M$ is a
matching in $K_n$. Since $|M|\geq k$, we can choose $M_1\subseteq M$
such that $|M_1|=k$. Let $M_1=\{u_iw_i|1\leq i\leq k\}$. Choose
$S=\{u_1,u_2,\cdots,u_k\}$. We will show that $\kappa(S)<
n-\frac{k}{2}-1$. Clearly, $|\bar{S}|=n-k$, and let
$\bar{S}=\{w_1,w_2,\cdots,w_{n-k}\}$. Since each tree in
$\mathscr{T}_2$ contains at least one vertex of $\bar{S}$, it
follows that $|\mathscr{T}_2|\leq n-k$. By the definition of
$\mathscr{T}_1$, we have $|\mathscr{T}_1|\leq \frac{k}{2}$. If
$|\mathscr{T}_1|\leq \frac{k}{2}-2$, then
$\kappa(S)\leq \lambda(S)=|\mathscr{T}|=|\mathscr{T}_1|+|\mathscr{T}_2|\leq
(\frac{k}{2}-2)+(n-k)=n-\frac{k}{2}-2<n-\frac{k}{2}-1$, as desired.
Let us assume $\frac{k}{2}-1\leq |\mathscr{T}_1|\leq \frac{k}{2}$.

Consider the case $|\mathscr{T}_1|=\frac{k}{2}-1$. Recall that
$|\mathscr{T}_2|\leq n-k$. Furthermore, we claim that
$|\mathscr{T}_2|\leq n-k-1$. Assume, to the contrary, that
$|\mathscr{T}_2|=n-k$. Let $T_1,T_2,\cdots,T_{n-k}$ be the $n-k$ edge-disjoint
$S$-Steiner trees in $\mathscr{T}_2$. For each tree $T_i \
(1\leq i\leq n-k)$, this tree only occupy one vertex of $\bar{S}$,
say $w_i$. Since $u_iw_i\in M_1 \ (1\leq i\leq k)$, namely,
$u_iw_i\notin E(G)$, and each $T_i\ (1\leq i\leq k)$ is an $S$-Steiner tree in $\mathscr{T}_2$, it follows that this tree
$T_i$ must contain at least one edge in $G[S]=K_k$. So the trees
$T_1,T_2,\cdots,T_{k}$ must use at least $k$ edges in $G[S]$, and
$|\mathscr{T}_1|=\frac{{{k}\choose{2}}-k}{k-1}=\frac{k-2}{2}-\frac{1}{k-1}$.
Since $|\mathscr{T}_1|$ is an integer, we have
$|\mathscr{T}_1|<\frac{k-2}{2}$, a contradiction. We conclude that
$|\mathscr{T}_2|\leq n-k-1$, and hence
$\kappa(S)\leq \lambda(S)=|\mathscr{T}|=|\mathscr{T}_1|+|\mathscr{T}_2|\leq
(\frac{k}{2}-1)+(n-k-1)=n-\frac{k}{2}-2<n-\frac{k}{2}-1$, as
desired.

Consider the case $|\mathscr{T}_1|=\frac{k}{2}$. We claim that
$|\mathscr{T}_2|\leq n-k-2$. Assume, to the contrary, that
$n-k-1\leq |\mathscr{T}_2|\leq n-k$. Since
$|\mathscr{T}_1|=\frac{k}{2}$, it follows that each edge of $G[S]$
is occupied by some tree in $\mathscr{T}_1$, which implies that each
tree in $\mathscr{T}_2$ only uses the edges of $E_G[S,\bar{S}]\cup
E(G[\bar{S}])$. Suppose that $T_1$ is a tree in $\mathscr{T}_2$
occupying $w_1$. Since $u_1w_1\notin E(G)$, if $T_1$ contains three
vertices of $\bar{S}$, then the remaining $n-k-3$ vertices in
$\bar{S}$ must be contained in at most $n-k-3$ trees in
$\mathscr{T}_2$, which results in $|\mathscr{T}_2|\leq
(n-k-3)+1=n-k-2$, a contradiction. So we assume that the tree $T_1$
contains another vertex of $\bar{S}$ except $w_1$, say $w_2$. Recall
that $k\geq 4$. Then $|\bar{S}|\geq k\geq 4$. By the same reason,
there is another tree $T_2$ containing two vertices of $\bar{S}$,
say $w_3,w_4$. Furthermore, the remaining $n-k-4$ vertices in
$\bar{S}$ must be contained in at most $n-k-4$ trees in
$\mathscr{T}_2$, which results in $|\mathscr{T}_2|\leq
(n-k-4)+2=n-k-2$, a contradiction. We conclude that
$|\mathscr{T}_2|\leq n-k-2$. Since $|\mathscr{T}_1|=\frac{k}{2}$, we
have $\kappa(S)\leq
\lambda(S)=|\mathscr{T}|=|\mathscr{T}_1|+|\mathscr{T}_2|\leq
\frac{k}{2}+(n-k-2)=n-\frac{k}{2}-2<n-\frac{k}{2}-1$, as
desired.\qed
\end{pf}

\begin{lem}\label{lem9}
If $n \ (n\geq 4)$ is even and $M$ is a set of edges in the complete
graph $K_n$ such that $1\leq |M|\leq n-1$ and $1\leq
\Delta(K_n[M])\leq \frac{n}{2}$, then $G=K_n- M$ contains
$\frac{n-2}{2}$ edge-disjoint spanning trees.
\end{lem}
\begin{pf}
Let $\mathscr{P}=\bigcup_{i=1}^pV_i$ be a partition of $V(G)$ with
$|V_i|=n_i \ (1\leq i\leq p)$, and $\mathcal {E}_p$ be the set of
edges between distinct blocks of $\mathscr{P}$ in $G$. It suffices
to show that $|\mathcal {E}_p|\geq \frac{n-2}{2}(|\mathscr{P}|-1)$
so that we can use Theorem \ref{th1}.

The case $p=1$ is trivial by Theorem \ref{th1}, thus we assume
$p\geq 2$. For $p=2$, we have $\mathscr{P}=V_1\cup V_2$. Set
$|V_1|=n_1$. Clearly, $|V_2|=n-n_1$. Since $\Delta(K_n[M])\leq
\frac{n}{2}$, it follows that $\delta(G)=n-1-\Delta(K_n[M])\geq
n-1-\frac{n}{2}=\frac{n-2}{2}$. Therefore, if $n_1=1$ then
$|\mathcal {E}_2|=|E_G[V_1,V_2]|\geq \frac{n-2}{2}$. Suppose
$n_1\geq 2$. Then $|\mathcal {E}_2|=|E_G[V_1,V_2]|\geq
{{n}\choose{2}}-(n-1)-{{n_1}\choose{2}}-{{n-n_1}\choose{2}}
=-n_1^2+nn_1-n+1$. Since $2\leq n_1 \leq n-2$, one can see that
$|\mathcal {E}_2|$ achives its minimum value when $n_1=2$ or
$n_1=n-2$. Thus $|\mathcal {E}_2|\geq n-3\geq \frac{n-2}{2}$ since
$n\geq 4$. The result follows from Theorem \ref{th1}.

Let us consider the remaining cases for $p$, namely, for $3\leq
p\leq n$. Since $|\mathcal {E}_p|\geq
{{n}\choose{2}}-|M|-\sum_{i=1}^p{{n_i}\choose{2}}\geq
{{n}\choose{2}}-(n-1)-\sum_{i=1}^p{{n_i}\choose{2}}
={{n-1}\choose{2}}-\sum_{i=1}^p{{n_i}\choose{2}}$, we only need to
show ${{n-1}\choose{2}}-\sum_{i=1}^p{{n_i}\choose{2}}\geq
\frac{n-2}{2}(p-1)$, that is, $(n-p)\frac{n-2}{2}\geq
\sum_{i=1}^p{{n_i}\choose{2}}$. Because
$\sum_{i=1}^p{{n_i}\choose{2}}$ achieves its maximum value when
$n_1=n_2=\cdots=n_{p-1}=1$ and $n_p=n-p+1$, we need inequality
$(n-p)\frac{n-2}{2}\geq {{1}\choose{2}}(p-1)+{{n-p+1}\choose{2}}$,
namely, $(n-p)\frac{p-3}{2}\geq 0$. It is easy to see that the
inequality holds since $3\leq p\leq n$. Thus, $|\mathcal {E}_p|\geq
{{n}\choose{2}}-|M|-\sum_{i=1}^p{{n_i}\choose{2}}\geq
\frac{n-2}{2}(p-1)$.

From Theorem $1$, there exist $\frac{n-2}{2}$ edge-disjoint spanning
trees in $G$, as desired.\qed
\end{pf}

\begin{lem}\label{lem10}
Let $k,n$ be two integers with $4\leq k\leq n$, and $M$ is an edge
set of the complete graph $K_n$ satisfying $\Delta(K_n[M])=1$. Then

$(1)$ If $|M|=k-1$, then $\kappa_k(K_n- M)\geq n-\frac{k}{2}-1$;

$(2)$ If $|M|=\lfloor \frac{n}{2} \rfloor$, then $\lambda_k(K_n-
M)\geq n-\frac{k}{2}-1$.
\end{lem}
\begin{pf}
$(1)$ Set $G=K_n- M$. Since $\Delta(K_n[M])=1$, it follows that $M$
is a matching of $K_n$. By the definition of $\kappa_k(G)$, we need
to show that $\kappa(S)\geq n-\frac{k}{2}-1$ for any $S\subseteq
V(G)$.

\textbf{Case 1}. There exists no $u,w$ in $S$ such that $uw\in M$.

Without loss of generality, let $S=\{u_1,u_2,\cdots, u_{k}\}$ such
that $u_1,u_2,\cdots,u_r$ are $M$-saturated but
$u_{r+1},u_{r+2},\cdots,u_k$ are $M$-unsaturated. Let
$M_1=\{u_iw_i\,|\,1\leq i\leq r\}\subseteq M$. Since $|M|=k-1$, it
follows that $0\leq r\leq k-1$. In this case, $u_iu_j\notin M \
(1\leq i,j\leq r)$. Clearly, $G[S]$ is a clique of order $k$. We
choose a path $P=u_{1}u_{2}\cdots u_{r}u_{r+1}$ in $G[S]$. Let
$G'=G- E(P)$. Then $G'[S]=K_k- E(P)$. Since $|E(P)|=r\leq k-1$ and
$\Delta(K_k[E(P)])=2\leq \frac{k}{2}$, it follows that $G'[S]$
contains $\frac{k-2}{2}$ edge-disjoint spanning trees, which are
also $\frac{k-2}{2}$ internally disjoint $S$-Steiner trees. These
trees together with the trees $T_i$ induced by the edges in
$\{u_1w_{i},u_2w_{i},
u_{i-1}w_{i},u_{i+1}w_{i},\cdots,u_{k}w_{i},u_{i}u_{i+1}\}  \ (1\leq
i\leq r)$ (see Figure 1 $(a)$) and the trees $T_j$ induced by the
edges in $\{u_1v_j,u_2v_j,\cdots, u_{k}v_j\}$ where $v_j\in \bar{S}-
\{w_1,w_2,\cdots,w_{r}\}=\{v_1,v_2,\cdots,v_{n-k-r}\}$ form
$\frac{k-2}{2}+r+(n-k-r)=n-\frac{k}{2}-1$ internally disjoint
$S$-Steiner trees. Thus, $\kappa(S)\geq n-\frac{k}{2}-1$, as
desired.

\textbf{Case 2}. There exist $u,w$ in $S$ such that $uw\in M$.

Without loss of generality, let $S=\{u_1,u_2,\cdots,
u_{r},u_{r+1},u_{r+2},\cdots,u_{r+s},u_{r+s+1},\cdots,$ $u_{k-r},
w_1, w_2,\cdots,w_{r}\}$ such that the vertices
$u_1,u_2,\cdots,u_{r+s}, w_1,w_2,\cdots,w_{r}$ are all $M$-saturated
and $u_iw_i\in M \ (1\leq i\leq r)$. Set $M_1=\{u_iw_i\,|\,1\leq
i\leq r\}$. In this case, $r\geq 1$ and $2r+s\leq k$. Since
$|M|=k-1$, it follows that $r+s\leq k-1$ and $s\leq k-2$.

First, we consider $2r+s=k$. Since $k$ is even, it follows that $s$
is even. If $s=0$, then $r=\frac{k}{2}$. Thus $S=\{u_1,u_2,\cdots,
u_{\frac{k}{2}},w_1,w_2,\cdots,w_{\frac{k}{2}}\}$. Clearly,
$M_1=\{u_iw_i\,|\,1\leq i\leq \frac{k}{2}\}$, $|M_1|=\frac{k}{2}\leq
k-1$ and $\Delta(K_n[M_1])=1<\frac{k}{2}$. By Lemma \ref{lem9},
$G[S]$ contains $\frac{k-2}{2}$ edge-disjoint spanning trees, which
are also $\frac{k-2}{2}$ internally disjoint $S$-Steiner trees. These trees together with the trees $T_j$ induced by
the edges in $\{u_1v_{j},u_2v_{j},\cdots,u_{\frac{k}{2}}v_{j}\}\cup
\{w_1v_{j},w_2v_{j},\cdots,w_{\frac{k}{2}}v_{j}\}$ form
$\frac{k-2}{2}+(n-k)$ internally disjoint $S$-Steiner trees, where $v_j\in \bar{S}=\{v_1,v_2,\cdots,v_{n-k}\}$. So,
$\kappa(S)\geq n-\frac{k}{2}-1$.

\begin{figure}[h,t,b,p]
\begin{center}
\scalebox{0.9}[0.9]{\includegraphics{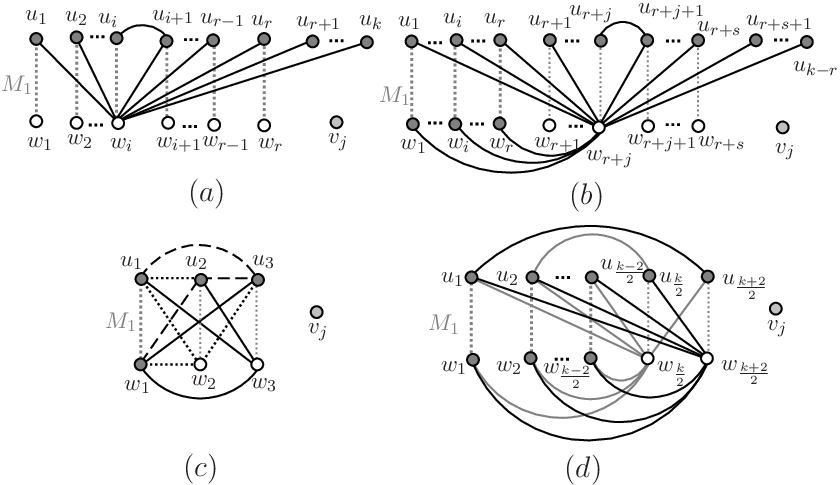}}\\[20pt]
 Figure 1. Graphs for $(1)$ of Lemma \ref{lem10}.
\end{center}\label{fig1}
\end{figure}

Consider $s=2$. Since $2r+s=k$, we have $r=\frac{k-2}{2}$. If $k=4$,
then $r=1$ and hence $S=\{u_1,u_2,u_3,w_1\}$. Clearly,
$M_1=\{u_1w_1\}$, and the tree $T_1$ induced by the edges in
$\{u_1u_2,u_1w_2,w_1w_2,u_3w_2\}$ and the tree $T_2$ induced by the
edges in $\{u_1u_3,u_2u_3,u_2w_1\}$ and the tree $T_3$ induced by
the edges in $\{u_1w_3,u_2w_3,w_1w_3,u_3w_1\}$ are three spanning
trees; see Figure 1 $(c)$. These trees together with the trees $T_j$
induced by the edges in $\{u_1v_{j},u_2v_{j}, u_{3}v_{j},w_1v_{j}\}$
form $3+(n-6)$ internally disjoint $S$-Steiner trees, where $v_j\in
\bar{S}- \{w_2,w_3\}=\{v_1,v_2,\cdots,v_{n-6}\}$. Thus,
$\kappa(S)\geq n-3=n-\frac{k}{2}-1$. Suppose $k\geq 6$. Then $r\geq
2$, $S=\{u_1,u_2,\cdots,
u_{\frac{k+2}{2}},w_1,w_2,\cdots,w_{\frac{k-2}{2}}\}$ and
$M_1=\{u_iw_i\,|\,1\leq i\leq \frac{k-2}{2}\}$. Clearly, the tree
$T_1$ induced by the edges in $\{u_1w_{\frac{k}{2}},$ $
u_2w_{\frac{k}{2}},\cdots,u_{\frac{k-2}{2}}w_{\frac{k}{2}},
u_{\frac{k+2}{2}}w_{\frac{k}{2}},u_2u_{\frac{k}{2}},
w_1w_{\frac{k}{2}},w_2w_{\frac{k}{2}},\cdots,$ $
w_{\frac{k-2}{2}}w_{\frac{k}{2}}\}$ and the tree $T_2$ induced by
the edges in $\{u_1w_{\frac{k+2}{2}}, u_2w_{\frac{k+2}{2}},\cdots,
u_{\frac{k}{2}}w_{\frac{k+2}{2}}\}\cup \{u_1u_{\frac{k+2}{2}},$
$w_1w_{\frac{k+2}{2}},
w_2w_{\frac{k+2}{2}},\cdots,w_{\frac{k-2}{2}}w_{\frac{k+2}{2}}\}$
are two internally disjoint $S$-Steiner trees; see Figure
1 $(d)$. Let $M_2=M_1\cup
\{u_1u_{\frac{k+2}{2}},u_2u_{\frac{k}{2}}\}$. Then
$|M_2|=|M_1|+2=\frac{k-2}{2}+2=\frac{k+2}{2}<k-1$ and
$\Delta(K_n[M_2])=2\leq \frac{k}{2}$, which implies that $G[S]-
\{u_1u_{\frac{k+2}{2}},u_2u_{\frac{k}{2}}\}=K_k- M_2$ contains
$\frac{k-2}{2}$ edge-disjoint spanning trees by Lemma \ref{lem9},
which are also $\frac{k-2}{2}$ internally disjoint $S$-Steiner trees. These trees together with $T_1,T_2$ and the trees
$T_j$ induced by the edges in $\{u_1v_j,u_2v_j,\cdots,
u_{\frac{k+2}{2}}v_j,w_1v_j,$ $w_2v_j,\cdots,u_{\frac{k-2}{2}}v_j\}$
are $\frac{k-2}{2}+2+(n-k-2)$ internally disjoint $S$-Steiner trees, where $v_j\in \bar{S}-
\{w_{\frac{k}{2}},w_{\frac{k+2}{2}}\}=\{v_1,v_2,\cdots,v_{n-k-2}\}$.
So, $\kappa(S)\geq n-\frac{k}{2}-1$.

Consider the remaining case for $s$, namely, for $4\leq s\leq k-2$.
Clearly, there exists a cycle of order $s$ containing
$u_{r+1},u_{r+2},\cdots,u_{r+s}$ in $K_k- M_1$, say
$C_{s}=u_{r+1}u_{r+2}\cdots u_{r+s}u_{r+1}$. Set $M'=M_1\cup
E(C_{s})$. Then $|M'|=r+s\leq k-1$ and $\Delta(K_n[M'])=2\leq
\frac{k}{2}$, which implies that $G- E(C_{s})=K_k- M'$ contains
$\frac{k-2}{2}$ edge-disjoint spanning trees by Lemma \ref{lem9}.
These trees together with the trees $T_{r+j}$ induced by the edges
in $\{u_1w_{r+j}, u_2w_{r+j},\cdots,$
$u_{r+j-1}w_{r+j},u_{r+j+1}w_{r+j},\cdots,
u_{r+s}w_{r+j},u_{r+j}u_{r+j+1},w_1w_{r+j},w_2w_{r+j},\cdots,
w_{r}w_{r+j}\} \ (1\leq j\leq s)$ form $\frac{k-2}{2}+s$ internally
disjoint trees; see Figure 2 $(b)$ (note that $u_{r+s}=u_{k-r}$).
These trees together with the trees $T_j'$ induced by the edges in
$\{u_1v_j,u_2v_j,\cdots,u_{r+s}v_j,w_{1}v_j,\cdots,w_{r}v_j\}$ form
$\frac{k-2}{2}+s+(n-2r-2s)=n-\frac{k}{2}-1$ internally disjoint
$S$-Steiner trees where $v_j\in \bar{S}-
\{w_{r+1},w_{r+2},\cdots,w_{r+s}\}=\{v_1,v_2,\cdots,v_{n-2r-2s}\}$.
Thus, $\kappa(S)\geq n-\frac{k}{2}-1$, as desired.

Next, assume $2r+s<k$. Then
$S=\{u_1,u_2,\cdots,u_{r+s},u_{r+s+1},\cdots,
u_{k-r},w_1,w_2,\cdots,w_{r}\}$ and $r+s+1\leq k-r$. If $s=0$, then
$S=\{u_1,u_2,\cdots, u_{k-r},w_1,w_2,\cdots,w_{r}\}$. Clearly,
$M_1=\{u_iw_i\,|\,1\leq i\leq r\}$, $|M_1|=r\leq k-1$ and
$\Delta(K_n[M_1])=1<\frac{k}{2}$. By Lemma \ref{lem9}, $G[S]$
contains $\frac{k-2}{2}$ edge-disjoint spanning trees. These trees
together with the trees $T_j$ induced by the edges in
$\{u_1v_{j},u_2v_{j},\cdots,u_{n-r}v_{j},w_1v_{j},w_2v_{j},\cdots,w_{r}v_{j}\}$
form $\frac{k-2}{2}+(n-k)$ internally disjoint $S$-Steiner trees, where
$v_j\in \bar{S}=\{v_1,v_2,\cdots,v_{n-k}\}$. Therefore,
$\kappa(S)\geq n-\frac{k}{2}-1$. Assume $s\geq 1$. Clearly, there
exists a path of length $s$ containing
$u_{r+1},u_{r+2},\cdots,u_{r+s},u_{r+s+1}$ in $G[S]$, say
$P_{s}=u_{r+1}u_{r+2}\cdots u_{r+s}u_{r+s+1}$. Set $M'=M_1\cup
E(P_{s})$. Then $|M'|=r+s\leq k-1$ and $\Delta(K_n[M'])=2\leq
\frac{k}{2}$, which implies that $G[S]-E(P_{s})=K_k- M'$ contains
$\frac{k-2}{2}$ edge-disjoint spanning trees by Lemma \ref{lem9},
which are also $\frac{k-2}{2}$ internally disjoint $S$-Steiner trees. These trees together with the trees $T_{r+j}$
induced by the edges in $\{u_1w_{r+j},
u_2w_{r+j},\cdots,u_{r+j-1}w_{r+j},u_{r+j+1}w_{r+j},\cdots,
u_{k-r}w_{r+j},u_{r+j}u_{r+j+1},w_1w_{r+j},w_2w_{r+j},$ $\cdots,
w_{r}w_{r+j}\} \ (1\leq j\leq s)$ form $\frac{k-2}{2}+s$ internally
disjoint $S$-Steiner trees; see Figure 1 $(b)$. These trees together
with the trees $T_j'$ induced by the edges in
$\{u_1v_j,u_2v_j,\cdots ,u_{k-r}v_j,$
$w_{1}v_j,w_{2}v_j,\cdots,w_{r}v_j\}$ form
$\frac{k-2}{2}+s+(n-k+r)-(r+s)=n-\frac{k}{2}-1$ internally disjoint
$S$-Steiner trees where $v_j\in \bar{S}-
\{w_{r+1},w_{r+2},\cdots,w_{r+s}\}=\{v_1,v_2,\cdots,v_{n-k-s}\}$.
So, $\kappa(S)\geq n-\frac{k}{2}-1$, as desired.

We conclude that $\kappa(S)\geq n-\frac{k}{2}-1$ for any $S\subseteq
V(G)$. From the arbitrariness of $S$, it follows that
$\kappa_k(G)\geq n-\frac{k}{2}-1$.

$(2)$ Set $G=K_n- M$. Assume that $n$ is even. Thus $M$ is a perfect
matching of $K_n$, and all vertices of $G$ are $M$-saturated. By the
definition of $\lambda_k(G)$, we need to show that $\lambda(S)\geq
n-\frac{k}{2}-1$ for any $S\subseteq V(G)$.

\textbf{Case 3}. There exists no $u,w$ in $S$ such that $uw\in M$.

Without loss of generality, let $S=\{u_1,u_2,\cdots, u_{k}\}$. In
this case, $u_iu_j\notin M \ (1\leq i,j\leq k)$. Let
$M_1=\{u_iw_i\,|\,1\leq i\leq k\}\subseteq M=\{u_iw_i\,|\,1\leq
i\leq \frac{n}{2}\}$. Clearly, $w_i\notin S \ (1\leq i\leq
\frac{n}{2})$ and $u_j\notin S \ (k+1\leq j\leq \frac{n}{2})$. Since
$G[S]$ is a clique of order $k$, it follows that there are
$\frac{k}{2}$ edge-disjoint spanning trees in $G[S]$, which are also
$\frac{k}{2}$ edge-disjoint $S$-Steiner trees. These
trees together with the trees $T_{i}$ induced by the edges in
$\{u_1w_{i},u_2w_{i},u_{i-1}w_{i},u_{i+1}w_{i},\cdots,
u_{k}w_{i},u_{i}w_{k},w_{i}w_{k}\} \ (1\leq i\leq k-1)$ (see Figure
2 $(a)$) and the trees $T_{j}'$ induced by the edges in
$\{u_1u_j,u_2u_j,\cdots,u_{k}u_j\} \ (k+1\leq j\leq \frac{n}{2})$
and the trees $T_{j}''$ induced by the edges in
$\{u_1w_j,u_2w_j,\cdots,u_{k}w_j\} \ (k+1\leq j\leq \frac{n}{2})$
form $\frac{k}{2}+(k-1)+(n-2k)=n-\frac{k}{2}-1$ edge-disjoint
$S$-Steiner trees. Therefore, $\lambda(S)\geq
n-\frac{k}{2}-1$, as desired.

\textbf{Case 4}. There exist $u,w$ in $S$ such that $uw\in M$.

Without loss of generality, let $S=\{u_1,u_2,\cdots,
u_{r+s},w_1,w_2,\cdots,w_{r}\}$ with $|S|=k=2r+s$, where $1\leq
r\leq \frac{k}{2}$ and $0\leq s\leq k-2$. Set
$M_1=\{u_iw_i\,|\,1\leq i\leq r\}\subseteq M=\{u_iw_i\,|\,1\leq
i\leq \frac{n}{2}\}$. We claim that $r+s\leq k-1$. Otherwise, let
$r+s=k$. Combining this with $2r+s=k$, we have $r=0$, a
contradiction. Since $k=2r+s$ and $k$ is even, it follows that $s$
is even.

\begin{figure}[h,t,b,p]
\begin{center}
\scalebox{0.9}[0.9]{\includegraphics{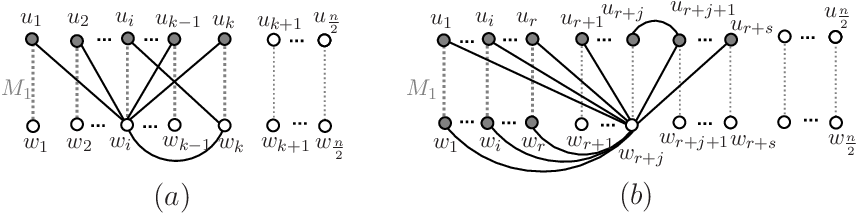}}\\[20pt]
 Figure 2. Graphs for $(2)$ of Lemma \ref{lem10}.
\end{center}\label{fig1}
\end{figure}

If $s=0$, then $r=\frac{k}{2}$. Clearly, $S=\{u_1,u_2,\cdots,
u_{\frac{k}{2}},w_1,w_2,\cdots,w_{\frac{k}{2}}\}$ and
$M_1=M=\{u_iw_i|1\leq i\leq \frac{k}{2}\}$. In addition, $|M_1|\leq
\frac{k}{2}<k-1$ and $\Delta(M\cap K_n[S])=1<\frac{k}{2}$. Then
$G[S]$ contains $\frac{k-2}{2}$ edge-disjoint spanning trees by
Lemma \ref{lem9}. These trees together with the trees $T_{i}$
induced by the edges in $\{u_1u_{i},u_2u_{i},
\cdots,u_{\frac{k}{2}}u_{i},w_1u_{i},w_2u_{i},\cdots,
w_{\frac{k}{2}}u_{i}\} \ (k+1\leq j\leq \frac{n}{2})$ and the trees
$T_{i}'$ induced by the edges in $\{u_1w_{i},u_2w_{i},
\cdots,u_{\frac{k}{2}}w_{i},w_1w_{i},w_2w_{i},$ $\cdots,
w_{\frac{k}{2}}w_{i}\} \ (\frac{k}{2}+1\leq i\leq \frac{n}{2})$ form
$n-\frac{k}{2}-1$ edge-disjoint $S$-Steiner trees. Thus,
$\lambda(S)\geq n-\frac{k}{2}-1$.

If $s=2$, then $r=\frac{k-2}{2}$. Then $S=\{u_1,u_2,\cdots,
u_{\frac{k+2}{2}},w_1,w_2,\cdots,w_{\frac{k-2}{2}}\}$ and
$M_1=\{u_iw_i\,|\,1\leq i\leq \frac{k-2}{2}\}\subseteq M$. If $k=4$,
then $r=1$ and hence $S=\{u_1,u_2,u_3,w_1\}$. Clearly,
$M_1=\{u_1w_1\}$, and the tree $T_{1}$ induced by the edges in
$\{u_1u_2,u_1w_2,w_1w_2,u_3w_2\}$ and the tree $T_{2}$ induced by
the edges in $\{u_1u_3,u_2u_3,u_2w_1\}$ and the tree $T_{3}$ induced
by the edges in $\{u_1w_3,u_2w_3,w_1w_3,u_3w_1\}$ are three
edge-disjoint spanning trees; see Figure 1 $(c)$. These trees
together with the trees $T_{j}$ induced by the edges in
$\{u_1u_{j},u_2u_{j}, u_{3}u_{j},w_1u_{j}\} \ (4\leq k\leq
\frac{n}{2})$ and the trees $T_{j}'$ induced by the edges in
$\{u_1w_{j},u_2w_{j}, u_{3}w_{j},w_1u_{j}\} \ (4\leq k\leq
\frac{n}{2})$ form $3+(n-6)$ edge-disjoint $S$-Steiner trees. So, $\lambda(S)\geq n-3=n-\frac{k}{2}-1$, as desired. Suppose
$k\geq 6$. Then $r\geq 2$, $S=\{u_1,u_2,\cdots,
u_{\frac{k+2}{2}},w_1,w_2,\cdots,w_{\frac{k-2}{2}}\}$ and
$M_1=\{u_iw_i|1\leq i\leq \frac{k-2}{2}\}$. Clearly, the tree
$T_{1}$ induced by the edges in $\{u_1w_{\frac{k}{2}},
u_2w_{\frac{k}{2}},\cdots,u_{\frac{k-2}{2}}w_{\frac{k}{2}},$ $
u_{\frac{k+2}{2}}w_{\frac{k}{2}},u_2u_{\frac{k}{2}},
w_1w_{\frac{k}{2}},w_2w_{\frac{k}{2}},\cdots,$ $
w_{\frac{k-2}{2}}w_{\frac{k}{2}}\}$ and the tree $T_{2}$ induced by
the edges in $\{u_1w_{\frac{k+2}{2}},$ $u_2w_{\frac{k+2}{2}},\cdots,
u_{\frac{k}{2}}w_{\frac{k+2}{2}},u_1u_{\frac{k+2}{2}},$ $
w_1w_{\frac{k+2}{2}},w_2w_{\frac{k+2}{2}},\cdots,
w_{\frac{k-2}{2}}w_{\frac{k+2}{2}}\}$ are two edge-disjoint $S$-Steiner
trees; see Figure 1 $(d)$. Let $M_2=M_1\cup
\{u_1u_{\frac{k+2}{2}},u_2u_{\frac{k}{2}}\}$. Then
$|M_2|=|M_1|+2=\frac{k-2}{2}+2=\frac{k+2}{2}<k-1$ and
$\Delta(K_n[M_2])=2\leq \frac{k}{2}$, which implies that $G[S]-
\{u_1u_{\frac{k+2}{2}},u_2u_{\frac{k}{2}}\}=K_k-M_2$ contains
$\frac{k-2}{2}$ edge-disjoint spanning trees by Lemma \ref{lem9}.
These trees together with $T_1,T_2$ and the trees $T_{j}$ induced by
the edges in $\{u_1u_j,u_2u_j,\cdots,u_{\frac{k+2}{2}}u_j,$
$w_1u_j,w_2u_j,\cdots,u_{\frac{k-2}{2}}u_j\} \ (\frac{k}{2}+2\leq
j\leq \frac{n}{2})$ and the trees $T_{j}'$ induced by the edges in
$\{u_1w_j,u_2w_j,\cdots,u_{\frac{k+2}{2}}w_j,w_1w_j,w_2w_j,\cdots,u_{\frac{k-2}{2}}w_j\}
\ (\frac{k}{2}+2\leq j\leq \frac{n}{2})$ are
$\frac{k-2}{2}+2+(n-k-2)$ edge-disjoint $S$-Steiner trees. Therefore, $\lambda(S)\geq n-\frac{k}{2}-1$, as desired.

Consider the remaining case $s$ with $4\leq s\leq k-2$. Clearly,
there exists a cycle of order $s$ containing
$u_{r+1},u_{r+2},\cdots,u_{r+s}$ in $K_k- M_1$, say
$C_{s}=u_{r+1}u_{r+2}\cdots u_{r+s}u_{r+1}$. Set $M'=M_1\cup
E(C_{s})$. Then $|M'|=r+s\leq k-1$ and $\Delta(K_n[M'])=2\leq
\frac{k}{2}$, which implies that $G- E(C_{s})$ contains
$\frac{k-2}{2}$ edge-disjoint spanning trees by Lemma \ref{lem9}.
These trees together with the trees $T_{r+j}$ induced by the edges
in $\{u_1w_{r+j},u_2w_{r+j},\cdots, u_{r+j-1}w_{r+j},$ $
u_{r+j+1}w_{r+j},\cdots,u_{r+s}w_{r+j},u_{r+j}u_{r+j+1}, w_1w_{r+j},
w_2w_{r+j},\cdots,w_{r}w_{r+j}\}  \ (1\leq j\leq s)$ form
$\frac{k-2}{2}+s$ edge-disjoint $S$-Steiner trees; see Figure 2 $(b)$.
These trees together with the trees $T_i'$ induced by the edges in
$\{u_1u_i,u_2u_i,\cdots,u_{r+s}u_i,w_{1}u_i,\cdots,w_{r}u_{i}\} \
(r+s+1\leq i\leq \frac{n}{2})$ and the trees $T_i''$ induced by the
edges in
$\{u_1w_i,u_2w_i,\cdots,u_{r+s}w_i,w_{1}w_i,\cdots,w_{r}w_{i}\} \
(r+s+1\leq i\leq \frac{n}{2})$ form
$(n-2r-2s)+(\frac{k-2}{2}+s)=n-\frac{k}{2}-1$ edge-disjoint $S$-Steiner
trees since $2r+s=k$. Thus, $\lambda(S)\geq
n-\frac{k}{2}-1$, as desired.

We conclude that $\lambda(S)\geq n-\frac{k}{2}-1$ for any
$S\subseteq V(G)$. From the arbitrariness of $S$, it follows that
$\lambda_k(G)\geq n-\frac{k}{2}-1$. For $n$ odd, $M$ is a maximum
matching and we can also check that $\lambda_k(G)\geq
n-\frac{k}{2}-1$ similarly.\qed
\end{pf}

\begin{lem}\label{lem11}
Let $n$ and $k$ be two integers such that $k$ is even and $4\leq
k\leq n$. If $M$ is a set of edges in the complete graph $K_n$ such
that $|M|=k-1$, and $2\leq \Delta(K_n[M])\leq \frac{k}{2}$, then
$\kappa_k(K_n- M)\geq n-\frac{k}{2}-1$.
\end{lem}
\begin{pf}
Set $G=K_n- M$. For $n=k$, there are $\frac{n-2}{2}$ edge-disjoint
spanning trees by Lemma \ref{lem9}, and hence $\kappa_n(G)=\lambda_n(G)\geq
\frac{n-2}{2}$. So from now on, we assume $n\geq k+1$. Let
$S=\{u_1,u_2,\cdots,u_k\}\subseteq V(G)$ and $\bar{S}=V(G)-
S=\{w_1,w_2,\cdots,w_{n-k}\}$. We have the following two cases to
consider.

\textbf{Case 1}. $M\subseteq E(K_n[S])\cup E(K_n[\bar{S}])$.

Let $M'=M\cap E(K_n[S])$ and $M''=M\cap E(K_n[\bar{S}])$. Then
$|M'|+|M''|=|M|=k-1$ and $0\leq |M'|,|M''|\leq k-1$. We can regard
$G[S]$ as a complete graph $K_k$ by deleting $|M'|$ edges. Since
$2\leq \Delta(K_n[M])\leq \frac{k}{2}$ and $M'\subseteq M$, it
follows that $\Delta(K_n[M'])\leq \Delta(K_n[M])\leq \frac{k}{2}$.
From Lemma \ref{lem9}, there exist $\frac{k-2}{2}$ edge-disjoint
spanning trees in $G[S]$. Actually, these $\frac{k-2}{2}$
edge-disjoint spanning trees are all internally disjoint $S$-Steiner
trees in $G[S]$. All these trees together with the trees $T_i$
induced by the edges in $\{w_iu_1,w_iu_2,\cdots, w_iu_{k}\} \ (1\leq
i\leq n-k)$ form $\frac{k-2}{2}+(n-k)=n-\frac{k}{2}-1$ internally
disjoint $S$-Steiner trees, and hence $\kappa(S)\geq
n-\frac{k}{2}-1$. From the arbitrariness of $S$, we have
$\kappa_k(G)\geq n-\frac{k}{2}-1$, as desired.

\textbf{Case 2}. $M\nsubseteq E(K_n[S])\cup E(K_n[\bar{S}])$.

In this case, there exist some edges of $M$ in $E_{K_n}[S,\bar{S}]$.
Let $M'=M\cap E(K_n[S])$, $M''=M\cap E(K_n[\bar{S}])$, and
$|M'|=m_1$ and $|M''|=m_2$. Clearly, $0\leq m_i\leq k-2 \ (i=1,2)$.
For $w_i\in \bar{S}$, let $|E_{K_n[M]}[w_i,S]|=x_i$, where $1\leq
i\leq n-k$. Without loss of generality, let $x_1\geq x_2\geq \cdots
\geq x_{n-k}$. Because there exist some edges of $M$ in
$E_{K_n}[S,\bar{S}]$, we have $x_1\geq 1$. Since $2\leq
\Delta(K_n[M])\leq \frac{k}{2}$, it follows that
$x_i=|E_{K_n[M]}[w_i,S]|\leq d_{K_n[M]}(w_i)\leq \Delta(K_n[M])\leq
\frac{k}{2}$ for $1\leq i\leq n-k$. We claim that there exists at
most one vertex in $K_n[M]$ such that its degree is $\frac{k}{2}$.
Assume, to the contrary, that there are two vertices, say $w$ and
$w'$, such that $d_{K_n[M]}(w)=d_{K_n[M]}(w')=\frac{k}{2}$. Then
$|M|\geq d_{K_n[M]}(w)+d_{K_n[M]}(w')=\frac{k}{2}+\frac{k}{2}=k$,
contradicting $|M|=k-1$. We conclude that there exists at most
one vertex in $K_n[M]$ such that its degree is $\frac{k}{2}$. Recall
that $x_{n-k}\leq x_{n-k-1}\leq \cdots \leq x_2\leq x_1\leq
\frac{k}{2}$. So $x_1=\frac{k}{2}$ and $x_i\leq \frac{k-2}{2} \
(2\leq i\leq n-k)$, or $x_i\leq \frac{k-2}{2} \ (1\leq i\leq n-k)$.
Since $|E_{K_n[M]}[w_i,S]|=x_i$, we have $|E_{G}[w_i,S]|=k-x_i$.
Since $2\leq \Delta(K_n[M])\leq \frac{k}{2}$, it follows that
$\delta(G[S])\geq k-1-\frac{k}{2}=\frac{k-2}{2}$.

Our basic idea is to seek for some edges in $G[S]$, and combine them
with the edges of $E_G[S,\bar{S}]$ to form $n-k$ internally disjoint
trees, say $T_1,T_2,\cdots,T_{n-k}$, with their roots
$w_1,w_2,\cdots,w_{n-k}$, respectively. Let $G'=G-
(\bigcup_{j=1}^{n-k}E(T_j))$. We will prove that $G'[S]$ satisfies
the conditions of Lemma \ref{lem9} so that $G'[S]$ contains
$\frac{k-2}{2}$ edge-disjoint spanning trees, which are also
$\frac{k-2}{2}$ internally disjoint $S$-Steiner trees.
These trees together with $T_1,T_2,\cdots,T_{n-k}$ are our
$n-\frac{k}{2}-1$ desired trees. Thus, $\kappa(S)\geq
n-\frac{k}{2}-1$. So we can complete our proof by the arbitrariness
of $S$.

For $w_1\in \bar{S}$, without loss of generality, let $S=S_1^1\cup
S_2^1$ and $S_1^1=\{u_1,u_2,\cdots,u_{x_1}\}$ such that $u_jw_1\in
M$ for $1\leq j\leq x_1$. Set $S_2^1=S-
S_1^1=\{u_{x_1+1},u_{x_1+2},\cdots,u_{k}\}$. Then $u_jw_1\in E(G)$
for $x_1+1\leq j\leq k$. One can see that the tree $T_1'$ induced by
the edges in $\{w_1u_{x_1+1}, w_1u_{x_1+2},\cdots,w_1u_{k}\}$ is a
Steiner tree connecting $S_2^1$. Our current idea is to seek for
$x_1$ edges in $E_G[S_1^1,S_2^1]$ and add them to $T_1'$ to form a
Steiner tree connecting $S$. For each $u_j\in S_1^1 \ (1\leq j\leq
x_1)$, we claim that $|E_G[u_j,S_2^1]|\geq 1$. Otherwise, let
$|E_G[u_j,S_2^1]|=0$. Then $|E_{K_n[M]}[u_j,S_2^1]|=k-x_1$ and hence
$|M|\geq |E_{K_n[M]}[u_j,S_2^1]|+d_{K_n[M]}(w_1)\geq (k-x_1)+x_1=k$,
which contradicts $|M|=k-1$. We conclude that for each $u_j\in S_1^1
\ (1\leq j\leq x_1)$ there is at least one edge in $G$ connecting it
to a vertex of $S_2^1$. Choose the vertex with the smallest
subscript among all the vertices of $S_1^1$ having maximum degree in
$G[S]$, say $u_{1}'$. Then we select the vertex adjacent to $u_1'$
with the smallest subscript among all the vertices of $S_2^1$ having
maximum degree in $G[S]$, say $u_1''$. Let $e_{11}=u_1'u_1''$.
Consider the graph $G_{11}=G- e_{11}$, and choose the vertex with
the smallest subscript among all the vertices of $S_1^1- u_1'$
having maximum degree in $G_{11}[S]$, say $u_{2}'$. Then we select
the vertex adjacent to $u_2'$ with the smallest subscript among all
the vertices of $S_2^1$ having maximum degree in $G_{11}[S]$, say
$u_2''$. Set $e_{12}=u_2'u_2''$. Consider the graph $G_{12}=G_{11}-
e_{12}=G- \{e_{11},e_{12}\}$. Choose the one with the smallest
subscript among all the vertices of $S_1^1- \{u_1',u_2'\}$ having
maximum degree in $G_{12}[S]$, say $u_{3}'$, and select the vertex
adjacent to $u_3'$ with the smallest subscript among all the
vertices of $S_2^1$ having maximum degree in $G_{12}[S]$, say
$u_3''$. Put $e_{13}=u_3'u_3''$. Consider the graph $G_{13}=G_{12}-
e_{11}=G- \{e_{11},e_{12},e_{13}\}$. For each $u_j\in S_1^1 \ (1\leq
j\leq x_1)$, we proceed to find $e_{14},e_{15},\cdots,e_{1x_1}$ in
the same way, and obtain graphs $G_{1j}=G-
\{e_{11},e_{12},\cdots,e_{1(j-1)}\} \ (1\leq j\leq x_1)$. Let
$M_1=\{e_{11},e_{12},\cdots,e_{1x_1}\}$ and $G_{1}=G- M_1$. Thus the
tree $T_1$ induced by the edges in
$\{w_1u_{x_2+1},w_1u_{x_2+2},\cdots,w_1u_{k}\}\cup \{e_{11},
e_{12},\cdots,e_{1x_1}\}$ is our desired tree.

Let us now prove the following claim.\vskip 0.5em

\textbf{Claim 1}. $\delta(G_{1}[S])\geq \frac{k-2}{2}$.\vskip 0.5em

\noindent{\itshape Proof of Claim $1$.} Assume, to the contrary,
that $\delta(G_{1}[S])\leq \frac{k-4}{2}$. Then there exists a
vertex $u_p\in S$ such that $d_{G_{1}[S]}(u_p)\leq \frac{k-4}{2}$.
If $u_p\in S_2^1$, then by our procedure
$d_{G[S]}(u_p)=d_{G_{1}[S]}(u_p)+1\leq \frac{k-2}{2}$, which implies
that $d_{M\cap K_n[S]}(u_p)\geq k-1-\frac{k-2}{2}=\frac{k}{2}$.
Since $w_1u_p\in M$, it follows that $d_{K_n[M]}(u_p)\geq d_{M\cap
K_n[S]}(u_p)+1\geq \frac{k+2}{2}$, which contradicts
$\Delta(K_n[M])\leq \frac{k}{2}$. Let us now assume $u_p\in S_2^1$.
By the above procedure, there exists a vertex $u_q\in S_1^1$ such
that when we select the edge $e_{1j}=u_pu_q \ (1\leq j\leq x_1)$
from $G_{1(j-1)}[S]$ the degree of $u_p$ in $G_{1j}[S]$ is equal to
$\frac{k-4}{2}$. Thus, $d_{G_{1j}[S]}(u_p)=\frac{k-4}{2}$ and
$d_{G_{1(j-1)}[S]}(u_p)=\frac{k-2}{2}$. From our procedure,
$|E_{G}[u_q,S_2^1]|=|E_{G_{1(j-1)}}[u_q,S_2^1]|$. Without loss of
generality, let $|E_{G}[u_q,S_2^1]|=t$ and $u_qu_j\in E(G)$ for
$x_1+1\leq j\leq x_1+t$; see Figure 3 $(a)$. Thus $u_p\in
\{u_{x_1+1},u_{x_1+2},\cdots,u_{x_1+t}\}$, and $u_qu_j\in M$ for
$x_1+t+1\leq j\leq k$. Because $|E_G[u_j,S_2^1]|\geq 1$ for each
$u_j\in S_1^1 \ (1\leq j\leq x_1)$, we have $t\geq 1$. Since
$|M|=k-1$ and $u_jw_1\in M$ for $1\leq j\leq x_1$, it follows that
$1\leq t\leq k-2$. Since $d_{G_{1(j-1)}[S]}(u_p)=\frac{k-2}{2}$, by
our procedure $d_{G_{1(j-1)}[S]}(u_j)\leq \frac{k-2}{2}$ for each
$u_j\in S_2^1 \ (x_1+1\leq j\leq x_1+t)$. Assume, to the contrary,
that there is a vertex $u_{s} \ (x_1+1\leq s\leq x_1+t)$ such that
$d_{G_{1(j-1)}[S]}(u_s)\geq \frac{k-2}{2}$. Then we should have
selected the edge $u_qu_s$ instead of $e_{1j}=u_qu_p$ by our
procedure, a contradiction. We conclude that
$d_{G_{1(j-1)}[S]}(u_r)\leq \frac{k-2}{2}$ for each $u_r\in S_1^1 \
(x_1+1\leq r\leq x_1+t)$. Clearly, there are at least
$k-1-\frac{k-2}{2}=\frac{k}{2}$ edges incident to each $u_r \
(x_1+1\leq r\leq x_1+t)$ belonging to $M\cup
\{e_{11},e_{12},\cdots,e_{1(j-1)}\}$. Since $j\leq x_1$ and
$u_qu_j\in M$ for $x_i+t+1\leq j\leq k$, we have
$$
|E_{K_n[M]}[u_q,S_2^1]|+\sum_{j=1}^td_{K_n[M]}(u_j)\geq
k-x_1-t+\frac{k}{2}t-(j-1)-{t\choose 2}=k+\frac{(k-2)}{2}t-x_1-j+1-{t\choose 2}
$$
and hence
\begin{eqnarray*}
|M|&\geq&|M\cap
(E_{K_n}[w_1,S])|+\sum_{j=1}^td_{K_n[M]}(u_j)+|E_{K_n[M]}[u_q,S_1^1]|\\
&\geq& x_1+\left(k+\frac{(k-2)}{2}t-x_1-j+1\right)-{t\choose 2}\\
&=&-\frac{t^2}{2}+\frac{t}{2}+\frac{(k-2)}{2}t+
k-j+1\\
&=&-\frac{t^2}{2}+\frac{(k-1)}{2}t+
k-j+1\\
&=&-\frac{1}{2}\left(t-\frac{k-1}{2}\right)^2+\frac{(k-1)^2}{8}+
k-j+1\\
&\geq&\frac{k}{2}-1+
k-j+1~~~~~~~~~~~~~~~~~~~~~~~~~~~~~~~(since~1\leq t\leq k-2)\\
&=&\frac{k}{2}+
k-j\\
&\geq&k,~~~~~~~~~~~~~~~~~~~~~~~~~~~~~~~~~~~~~~~~~~~~~\left(since~j\leq x_1~and~x_1\leq \frac{k}{2}\right)\\
\end{eqnarray*}
contradicting $|M|=k-1$.

\begin{figure}[h,t,b,p]
\begin{center}
\scalebox{0.9}[0.9]{\includegraphics{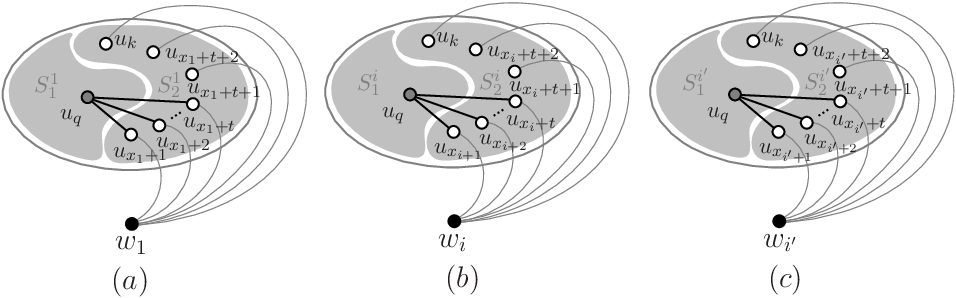}}\\
Figure 3. Graphs for Lemma \ref{lem11}.
\end{center}
\end{figure}

By Claim $1$, we have $\delta(G_{1}[S])\geq \frac{k-2}{2}$. Recall
that there exists at most one vertex in $K_n[M]$ such that its
degree is $\frac{k}{2}$, and $x_{n-k}\leq x_{n-k-1}\leq \cdots \leq
x_2\leq x_1\leq \frac{k}{2}$. Then $x_i\leq \frac{k-2}{2}$ for
$2\leq i\leq n-k$. Now we continue to introduce our procedure.

For $w_2\in \bar{S}$, without loss of generality, let $S=S_1^2\cup
S_2^2$ and $S_1^2=\{u_1,u_2,\cdots,u_{x_2}\}$ such that $u_jw_2\in
M$ for $1\leq j\leq x_2$. Let $S_2^2=S-
S_1^2=\{u_{x_2+1},u_{x_2+2},\cdots,u_{k}\}$. Then $u_jw_2\in E(G)$
for $x_2+1\leq j\leq k$. Clearly, the tree $T_2'$ induced by the
edges in $\{w_2u_{x_2+1}, w_2u_{x_2+2},\cdots,w_2u_{k}\}$ is a
Steiner tree connecting $S_2^2$. Our idea is to seek for $x_2$ edges
in $E_{G_1}[S_1^2,S_2^2]$ and add them to $T_2'$ to form a Steiner
tree connecting $S$. For each $u_j\in S_1^2 \ (1\leq j\leq x_2)$, we
claim that $|E_{G_1}[u_j,S_2^2]|\geq 1$. Otherwise, let
$|E_{G_1}[u_j,S_2^2]|=0$. Recall that $|M_1|=x_1$. Then there exist
$k-x_2$ edges between $u_j$ and $S_2^2$ belonging to $M\cup M_1$,
and hence $|E_{K_n[M]}[u_j,S_2^2]|\geq k-x_2-x_1$. Therefore,
$|M|\geq |E_{K_n[M]}[u_j,S_2^2]|+d_{K_n[M]}(w_1)+d_{K_n[M]}(w_2)
\geq (k-x_2-x_1)+x_1+x_2=k$, which contradicts $|M|=k-1$. Choose the
vertex with the smallest subscript among all the vertices of $S_1^2$
having maximum degree in $G_1[S]$, say $u_{1}'$. Then we select the
vertex adjacent to $u_1'$ with the smallest subscript among all the
vertices of $S_2^2$ having maximum degree in $G_1[S]$, say $u_1''$.
Let $e_{21}=u_1'u_1''$. Consider the graph $G_{21}=G_1- e_{21}$, and
choose the one with the smallest subscript among all the vertices of
$S_1^2- u_1'$ having maximum degree in $G_{21}[S]$, say $u_{2}'$.
Then we select the vertex adjacent to $u_2'$ with the smallest
subscript among all the vertices of $S_2^2$ having maximum degree in
$G_{21}[S]$, say $u_2''$. Set $e_{22}=u_2'u_2''$. Consider the graph
$G_{22}=G_{21}- e_{22}=G_1- \{e_{21},e_{22}\}$. For each $u_j\in
S_1^2 \ (1\leq j\leq x_2)$, we proceed to find
$e_{23},e_{24},\cdots,e_{2x_2}$ in the same way, and get graphs
$G_{2j}=G_1- \{e_{21},e_{22},\cdots,e_{2(j-1)}\} \ (1\leq j\leq
x_2)$. Let $M_2=\{e_{21},e_{22},\cdots,e_{2x_2}\}$ and $G_2=G_1-
M_1$. Thus the tree $T_2$ induced by the edges in
$\{w_2u_{x_2+1},w_2u_{x_2+2},\cdots,w_2u_{k}\}\cup \{e_{21},
e_{22},\cdots,e_{2x_2}\}$ is our desired tree. Furthermore, $T_2$
and $T_1$ are two internally disjoint $S$-Steiner trees.

For $w_i\in \bar{S}$, without loss of generality, let $S=S_1^i\cup
S_2^i$ and $S_1^i=\{u_1,u_2,\cdots,u_{x_i}\}$ such that $u_jw_i\in
M$ for $1\leq j\leq x_i$. Set $S_2^i=S-
S_1^i=\{u_{x_i+1},u_{x_i+2},\cdots,u_{k}\}$. Then $u_jw_i\in E(G)$
for $x_i+1\leq j\leq k$. One can see that the tree $T_i'$ induced by
the edges in $\{w_iu_{x_i+1}, w_iu_{x_i+2},\cdots,w_iu_{k}\}$ is a
Steiner tree connecting $S_2^i$. Our idea is to seek for $x_i$ edges
in $E_{G_{i-1}}[S_1^2,S_2^2]$ and add them to $T_i'$ to form a
Steiner tree connecting $S$. For each $u_j\in S_1^i \ (1\leq j\leq
x_i)$, we claim that $|E_{G_{i-1}}[u_j,S_2^i]|\geq 1$. Otherwise,
let $|E_{G_{i-1}}[u_j,S_2^i]|=0$. Recall that $|M_j|=x_j \ (1\leq
j\leq i)$. Then there are $k-x_i$ edges between $u_j$ and $S_2^i$
belonging to $M\cup (\bigcup_{j=1}^{i-1}M_j)$, and hence
$|E_{K_n[M]}[u_j,S_2^i]|\geq k-x_i-\sum_{j=1}^{i-1}x_j$. Therefore,
$|M|\geq |E_{K_n[M]}[u_j,S_2^i]|+\sum_{j=1}^{i}|M\cap
(K_n[w_j,S])|\geq k-x_i-\sum_{j=1}^{i-1}x_j+\sum_{j=1}^{i}x_j=k$,
contradicting $|M|=k-1$. Choose the vertex with the smallest
subscript among all the vertices of $S_1^i$ having maximum degree in
$G_{i-1}[S]$, say $u_{1}'$. Then we select the vertex adjacent to
$u_1'$ with the smallest subscript among all the vertices of $S_2^i$
having maximum degree in $G_{i-1}[S]$, say $u_1''$. Let
$e_{i1}=u_1'u_1''$. Consider the graph $G_{i1}=G_{i-1}- e_{i1}$,
choose the vertex with the smallest subscript among all the vertices
of $S_1^i- u_1'$ having maximum degree in $G_{i1}[S]$, say $u_{2}'$.
Then we select the vertex adjacent to $u_2'$ with the smallest
subscript among all the vertices of $S_2^i$ having maximum degree in
$G_{i1}[S]$, say $u_2''$. Set $e_{i2}=u_2'u_2''$. Consider the graph
$G_{i2}=G_{i1}- e_{i2}=G_{i-1}- \{e_{i1},e_{i2}\}$. For each $u_j\in
S_1^i \ (1\leq j\leq x_i)$, we proceed to find
$e_{i3},e_{i4},\cdots,e_{ix_i}$ in the same way, and get graphs
$G_{ij}=G_{i-1}- \{e_{i1},e_{i2},\cdots,e_{i(j-1)}\} \ (1\leq j\leq
x_i)$. Let $M_i=\{e_{i1},e_{i2},\cdots,e_{ix_2}\}$ and $G_i=G_{i-1}-
M_i$. Thus the tree $T_i$ induced by the edges in
$\{w_iu_{x_2+1},w_iu_{x_2+2},\cdots,w_iu_{k}\}\cup \{e_{i1},
e_{i2},\cdots,e_{ix_i}\}$ is our desired tree. Furthermore,
$T_1,T_2,\cdots,T_{i}$ are pairwise internally disjoint $S$-Steiner
trees.

We continue this procedure until we obtain $n-k$ pairwise internally
disjoint trees $T_1,T_2,\cdots, T_{n-k}$. Note that if there exists
some $x_j$ such that $x_j=0$ then $x_{j+1}=x_{j+2}=\cdots=x_{n-k}=0$
since $x_1\geq x_2\geq \cdots \geq x_{n-k}$. Then the trees $T_i$
induced by the edges in $\{w_iu_1,w_iu_2,\cdots,w_iu_{k}\} \ (j\leq
i\leq n-k)$ is our desired tree. From the above procedure, the
resulting graph must be $G_{n-k}=G- \bigcup_{i=1}^{n-k}M_i$. Let us
show the following claim.\vskip 0.5em

\textbf{Claim 2}. $\delta(G_{n-k}[S])\geq \frac{k-2}{2}$.\vskip 0.5em

\noindent{\itshape Proof of Claim $2$.} Assume, to the contrary,
that $\delta(G_{n-k}[S])\leq \frac{k-4}{2}$, namely, there exists a
vertex $u_p\in S$ such that $d_{G_{n-k}[S]}(u_p)\leq \frac{k-4}{2}$.
Since $\delta(G[S])\geq \frac{k-2}{2}$, by our procedure there
exists an edge $e_{ij}$ in $G_{i(j-1)}$ incident to the vertex $u_p$
such that when we pick up this edge,
$d_{G_{ij}[S]}(u_p)=\frac{k-4}{2}$ but
$d_{G_{i(j-1)}[S]}(u_p)=\frac{k-2}{2}$.

First, we consider the case $u_p\in S_2^i$. Then there exists a
vertex $u_q\in S_1^i$ such that when we select the edge
$e_{ij}=u_pu_q$ from $G_{i(j-1)}[S]$ the degree of $u_p$ in
$G_{ij}[S]$ is equal to $\frac{k-4}{2}$. Thus,
$d_{G_{ij}[S]}(u_p)=\frac{k-4}{2}$ and
$d_{G_{i(j-1)}[S]}(u_p)=\frac{k-2}{2}$. From our procedure,
$|E_{G_{i-1}}[u_q,S_2^i]|=|E_{G_{i(j-1)}}[u_q,S_2^i]|$. Without loss
of generality, let $|E_{G_{i-1}}[u_q,S_2^i]|=t$ and $u_qu_j\in
E(G_{i-1})$ for $x_i+1\leq j\leq x_i+t$; see Figure 3 $(b)$. Thus
$u_p\in \{u_{x_i+1},u_{x_i+2},\cdots,u_{x_i+t}\}$, and $u_qu_j\in
M\cup (\bigcup_{r=1}^{i-1}M_r)$ for $x_i+t+1\leq j\leq k$. Since
$x_i\leq \frac{k-2}{2} \ (2\leq i\leq n-k)$, it follows that
$|S_1^i|\leq \frac{k-2}{2}$. From this together with
$\delta(G_{i-1}[S])\geq \frac{k-2}{2}$, we have
$|E_{G_{i-1}}[u_q,S_1^i]|\geq 1$, that is, $t\geq 1$. Since
$d_{G_{i(j-1)}[S]}(u_p)=\frac{k-2}{2}$, by our procedure
$d_{G_{i(j-1)}[S]}(u_j)\leq \frac{k-2}{2}$ for each $u_j\in S_2^i \
(x_i+1\leq j\leq x_i+t)$. Assume, to the contrary, that there exists
a vertex $u_{s} \ (x_i+1\leq s\leq x_i+t)$ such that
$d_{G_{i(j-1)}[S]}(u_s)\geq \frac{k-2}{2}$. Then we should have
selected the edge $u_qu_s$ instead of $e_{ij}=u_qu_p$ by our
procedure, a contradiction. We conclude that
$d_{G_{i(j-1)}[S]}(u_r)\leq \frac{k-2}{2}$ for each $u_r\in S_2^i \
(x_i+1\leq r\leq x_i+t)$. Clearly, there are at least
$k-1-\frac{k-2}{2}=\frac{k}{2}$ edges incident to each $u_r \
(x_i+1\leq r\leq x_i+t)$ belonging to $M\cup
(\bigcup_{j=1}^{i-1}M_j) \bigcup
\{e_{i1},e_{i2},\cdots,e_{i(j-1)}\}$. Since $j\leq x_i$ and
$u_qu_j\in M\cup (\bigcup_{r=1}^{i-1}M_r)$ for $x_i+t+1\leq j\leq
k$, we have
\begin{eqnarray*}
&&|E_{K_n[M]}[u_q,S_2^i]|+\sum_{j=1}^td_{K_n[M]}(u_j)\\
&\geq&k-x_i-t+\frac{k}{2}t-\sum_{j=1}^{i-1}x_j-(j-1)-{t\choose 2}\\
&\geq&k+\frac{(k-2)}{2}t-\sum_{j=1}^{i}x_j-x_i+1-{t\choose 2}~~~~~~~~~~~~~~~~~~~~~~~~~(since~j\leq
x_i)\\
&=&-\frac{t^2}{2}+\frac{(k-1)}{2}t+k-\sum_{j=1}^{i}x_j-x_i+1\\
&=&-\frac{1}{2}\left(t-\frac{k-1}{2}\right)^2+\frac{(k-1)^2}{8}+k-\sum_{j=1}^{i}x_j-x_i+1\\
\end{eqnarray*}
and hence
\begin{eqnarray*}
|M|&\geq &\sum_{j=1}^{i}|M\cap
(E_{K_n}[w_j,S])|+\sum_{j=1}^td_{K_n[M]}(u_j)+|E_{K_n[M]}[u_q,S_2^i]|\\
&\geq &\sum_{j=1}^{i}x_j-\frac{1}{2}\left(t-\frac{k-1}{2}\right)^2+\frac{(k-1)^2}{8}+k-\sum_{j=1}^{i}x_j-x_i+1\\
&=&-\frac{1}{2}\left(t-\frac{k-1}{2}\right)^2+\frac{(k-1)^2}{8}+k-x_i+1\\
&\geq&\frac{k}{2}-1+k-x_i+1~~~~~~~~~~~~~~~~~~~~~~~~~(since~1\leq t\leq k-2)\\
&\geq&\frac{k}{2}+k-x_i\\
&\geq &
k+1,~~~~~~~~~~~~~~~~~~~~~~~~~~~~~~~~~~~~~~~~\left(since~x_i\leq
\frac{k-2}{2}\right)
\end{eqnarray*}
which contradicts $|M|=k-1$.

Next, assume $u_p\in S_1^i$. Recall that
$d_{G_{ij}[S]}(u_p)=\frac{k-4}{2}$. Since $u_p\in S_1^i$, it follows
that $d_{G_{i-1}[S]}(u_p)=\frac{k-2}{2}$. If $u_p\in \bigcap_{j=1}^i
S_1^j$, namely, $u_pw_j\in M \ (1\leq j\leq i)$, then by our
procedure $d_{G[S]}(u_p)=\frac{k-2}{2}+i-1$ and hence $d_{K_n[S]\cap
M}(u_p)=k-1-(\frac{k-2}{2}+i-1)=\frac{k}{2}-i+1$. Since $u_pw_j\in
M$ for each $w_j\in \bar{S} \ (1\leq j\leq i)$, we have
$d_{K_n[M]}(u_p)=d_{K_n[S]\cap M}(u_p)+d_{K_n[S,\bar{S}]\cap
M}(u_p)\geq (\frac{k}{2}-i+1)+i=\frac{k+2}{2}$, contradicting 
$\Delta(K_n[M])\leq \frac{k}{2}$. Combining this with $u_p\in
S_1^i$, we have $u_p\notin \bigcap_{j=1}^{i-1} S_1^i$ and we can
assume that there exists an integer $i' \ (i'\leq i-1)$ satisfying
the following conditions:

$\bullet$ $u_p\in S_2^{i'}$ and
$d_{G_{i'}[S]}(u_p)<d_{G_{i'-1}[S]}(u_p)$;

$\bullet$ if $u_p$ belongs to some $S_2^j \ (i'+1\leq j\leq i)$ then
$d_{G_j[S]}(u_p)=d_{G_{j-1}[S]}(u_p)$.

The above two conditions can be restated as follows:

$\bullet$ $u_pw_{i'}\in E(G)$ and
$d_{G_{i'}[S]}(u_p)<d_{G_{i'-1}[S]}(u_p)$;

$\bullet$ if $u_pw_j\in E(G) \ (i'+1\leq j\leq i)$ then
$d_{G_j[S]}(u_p)=d_{G_{j-1}[S]}(u_p)$.

In fact, we can find the integer $i'$ such that $u_pw_{i'}\in E(G)$
and $d_{G_{i'}[S]}(u_p)<d_{G_{i'-1}[S]}(u_p)$. Assume, to the
contrary, that for each $w_j \ (1\leq j\leq i)$, $u_pw_j\in M$, or
$u_pw_j\in E(G)$ but $d_{G_j[S]}(u_p)=d_{G_{j-1}[S]}(u_p)$. Let $i_1
\ (i_1\leq i)$ be the number of vertices nonadjacent to $u_p\in S$
in $\{w_1,w_2,\cdots,w_{i-1}\}\subseteq \bar{S}$. Without loss of
generality, let $w_{j}u_p\in M \ (1\leq j\leq i_1)$. Recall that
$d_{G_{ij}[S]}(u_p)=\frac{k-4}{2}$. Thus
$d_{G[S]}(u_p)=\frac{k-4}{2}+i_1$ and hence $d_{K_n[S]\cap
M}(u_p)\geq k-1-(\frac{k-4}{2}+i_1)=\frac{k+2}{2}-i_1$. Since
$w_ju_p\in M \ (1\leq j\leq i_1)$, it follows that $d_{K_n[S,\bar
S]\cap M}(u_p)\geq i_1$, which results in
$d_{K_n[M]}(u_p)=d_{K_n[S]\cap M}(u_p)+d_{K_n[S,\bar{S}]\cap
M}(u_p)\geq (\frac{k+2}{2}-i_1)+i_1=\frac{k+2}{2}$, contradicting 
$\Delta(K_n[M])\leq \frac{k}{2}$.

Now we turn our attention to $u_p\in S_2^{i'}$. Without loss of
generality, let $u_pw_j\in M \ (j\in
\{j_{1},j_{2},\cdots,j_{i_1}\})$, namely, $u_p\in S_1^{j_1}\cap
S_1^{j_2}\cap \cdots \cap S_1^{j_{i_1}}$, where
$j_{1},j_{2},\cdots,j_{i_1}\in \{i'+1,i'+2,\cdots,i\}$. Then
$u_pw_j\in E(G) \ (j\in \{i'+1,i'+2,\cdots,i\}-
\{j_{1},j_{2},\cdots,j_{i_1}\})$. Clearly, $i_1\leq i-i'$. Recall
that $u_p\in S_1^i$ and $d_{G_{ij}[S]}(u_p)=\frac{k-4}{2}$. Thus
$d_{G_{i'}[S]}(u_p)=\frac{k-4}{2}+i_1$. By our procedure, there
exists a vertex $u_q\in S_1^{i'}$ such that when we select the edge
$e_{i'j}=u_pu_q$ from $G_{i'(j-1)}[S]$ the degree of $u_p$ in
$G_{i'j}[S]$ is equal to $\frac{k-4}{2}+i_1$, that is,
$d_{G_{i'j}[S]}(u_p)=\frac{k-4}{2}+i_1$ and
$d_{G_{i'(j-1)}[S]}(u_p)=\frac{k-2}{2}+i_1$. From our procedure,
$|E_{G_{i'-1}}[u_q,S_2^{i'}]|=|E_{G_{i'(j-1)}}[u_q,S_2^{i'}]|$.
Without loss of generality, let $|E_{G_{i'-1}}[u_q,S_2^{i'}]|=t$ and
$u_qu_j\in E(G_{i'-1})$ for $x_{i'}+1\leq j\leq x_{i'}+t$; see
Figure 3 $(c)$. Thus $u_p\in
\{u_{x_{i'}+1},u_{x_{i'}+2},\cdots,u_{x_{i'}+t}\}$, and $u_qu_j\in
M\cup (\bigcup_{r=1}^{i'-1}M_r)$ for $x_{i'}+t+1\leq j\leq k$. Since
$x_{j}\leq \frac{k-2}{2} \ (2\leq j\leq n-k)$, it follows that
$|S_1^{i'}|\leq \frac{k-2}{2}$. From this together with
$\delta(G_{i'-1}[S])\geq \frac{k-2}{2}$, we have
$|E_{G_{i'-1}}[u_q,S_1^{i'}]|\geq 1$, that is, $t\geq 1$. Since
$d_{G_{i'(j-1)}[S]}(u_p)=\frac{k-2}{2}+i_1$, by our procedure
$d_{G_{i'(j-1)}[S]}(u_j)\leq \frac{k-2}{2}+i_1$ for each $u_j\in
S_2^{i'} \ (x_{i'}+1\leq j\leq x_{i'}+t)$. Assume, to the contrary,
that there is a vertex $u_{s} \ (x_{i'}+1\leq s\leq x_{i'}+t)$ such
that $d_{G_{i'(j-1)}[S]}(u_s)\geq \frac{k-2}{2}+i_1+1$. Then we
should have selected the edge $u_qu_s$ instead of $e_{i'j}=u_qu_p$
by our procedure, a contradiction. We conclude that
$d_{G_{i'(j-1)}[S]}(u_r)\leq \frac{k-2}{2}+i_1$ for each $u_r\in
S_2^{i'} \ (x_{i'}+1\leq r\leq x_{i'}+t)$. Clearly, there are at
least $k-1-(\frac{k-2}{2}+i_1)=\frac{k}{2}-i_1$ edges incident to
each $u_r \ (x_{i'}+1\leq r\leq x_{i'}+t)$ belonging to $M\cup
(\bigcup_{j=1}^{i'-1}M_j) \bigcup
\{e_{i'1},e_{i'2},\cdots,e_{i'(j-1)}\}$. Since $j\leq x_{i'}$ and
$u_qu_j\in M\cup (\bigcup_{r=1}^{i'-1}M_r)$ for $x_{i'}+t+1\leq
j\leq k$, we have
\begin{eqnarray*}
& &|E_{K_n[M]}[u_q,S_2^{i'}]|+\sum_{j=1}^td_{K_n[M]}(u_j)\\
&\geq& k-x_{i'}-t+\left(\frac{k}{2}-i_1\right)t-\sum_{j=1}^{i'-1}x_j-(j-1)-{t\choose 2}\\
&\geq& k-\sum_{j=1}^{i'}x_j+\left(\frac{k-2}{2}-i_1\right)t-x_{i'}+1-\frac{t(t-1)}{2}~~~~~~~~(since~j\leq x_{i'})\\
&=&-\frac{t^2}{2}+\frac{t}{2}+
k-\sum_{j=1}^{i'}x_j+\left(\frac{k-2}{2}-i+i'\right)t-x_{i'}+1~~~~(since~i_1\leq
i-i')\\
&=&-\frac{t^2}{2}+\left(\frac{k-1}{2}-i+i'\right)t+
k-\sum_{j=1}^{i'}x_j-x_{i'}+1\\
&=&-\frac{1}{2}\Big(t^2-(k-1-2i+2i')t\Big)+
k-\sum_{j=1}^{i'}x_j-x_{i'}+1\\
&=&-\frac{1}{2}\left(t-\frac{k-1-2i+2i'}{2}\right)^2+\frac{(k-1-2i+2i')^2}{8}+
k-\sum_{j=1}^{i'}x_j-x_{i'}+1\\
\end{eqnarray*}
and hence
\begin{eqnarray*}
|M|&\geq &\sum_{j=1}^{i}|M\cap
(E_{K_n}[w_j,S])|+\sum_{j=1}^pd_{K_n[M]}(u_j)+|E_{K_n[M]}[u_q,S_2^i]|\\
&\geq &\sum_{j=1}^{i}x_j-\frac{1}{2}\left(t-\frac{k-1-2i+2i'}{2}\right)^2+\frac{(k-1-2i+2i')^2}{8}+
k-\sum_{j=1}^{i'}x_j-x_{i'}+1\\
&=&-\frac{1}{2}\left(t-\frac{k-1-2i+2i'}{2}\right)^2+\frac{(k-1-2i+2i')^2}{8}+k+\sum_{j=i'+1}^{i}x_j-x_{i'}+1\\
&\geq&\frac{k}{2}-1-i+i'+k+\sum_{j=i'+1}^{i}x_j-x_{i'}+1~~~~~~~~~~~~~~~~~~(since~1\leq t\leq k-2~and\\
&&~~~~~~~~~~~~~~~~~~~~~~~~~~~~~~~~~~~~~~~~~~~~~~~~~~~~~~~~~~~~~~~~~~~~~k-1-2i+2i'\leq k-2)\\
&\geq & k,~~~~~~~~~~~~~~~\left(since~x_{i'}\leq
\frac{k-2}{2}~and~x_j\geq 1~for~i'+1\leq j\leq i\right)
\end{eqnarray*}
contradicting $|M|=k-1$. This completes the proof of Claim $2$.

From our procedure, we get $n-k$ internally disjoint Steiner trees
connecting $S$ in $G$, say $T_1,T_2,\cdots,T_{n-k}$. Recall that
$G_{n-k}=G-(\bigcup_{i=1}^{n-k}M_i)$. We can regard
$G_{n-k}[S]=G[S]- (\bigcup_{i=1}^{n-k}M_i)$ as a graph obtained from
the complete graph $K_k$ by deleting $|M'|+\sum_{i=1}^{n-k}|M_i|$
edges. Since
$|M'|+\sum_{i=1}^{n-k}|M_i|+|M''|=m_1+\sum_{i=1}^{n-k}x_i+m_2=k-1$,
we have $1\leq \sum_{i=1}^{n-k}|M_i|+m_1\leq k-1$. By Claim $2$,
$\delta(G_{n-k}[S])\geq \frac{k-2}{2}$ and hence $2\leq
\Delta(\overline{G_{n-k}}[S])\leq \frac{k}{2}$. From Lemma
\ref{lem9}, there exist $\frac{k-2}{2}$ edge-disjoint spanning trees
connecting $S$ in $G_{n-k}[S]$. These trees together with
$T_1,T_2,\cdots,T_{n-k}$ are $n-\frac{k}{2}-1$ internally disjoint
Steiner trees connecting $S$ in $G$. Thus, $\kappa(S)\geq
n-\frac{k}{2}-1$. From the arbitrariness of $S$, we have
$\kappa_k(G)\geq n-\frac{k}{2}-1$, as desired.\qed
\end{pf}\\

We are now in a position to prove our main results.\vskip 0.5em

\noindent \textbf{Proof of Theorem \ref{th2}}. Assume that
$\kappa_k(G)=n-\frac{k}{2}-1$. Since $G$ of order $n$ is connected,
we can regard $G$ as a graph obtained from the complete graph $K_n$
by deleting some edges. From Lemma \ref{lem4}, it follows that
$|M|\geq 1$ and hence $\Delta(K_n[M])\geq 1$. If $G=K_n-M$ where
$M\subseteq E(K_n)$ such that $\Delta(K_n[M])\geq \frac{k}{2}+1$,
then $\kappa_k(G)\leq \lambda_k(G)<n-\frac{k}{2}-1$ by Observation
\ref{obs1} and Corollary \ref{cor1}, a contradiction. So $1\leq
\Delta(K_n[M])\leq \frac{k}{2}$. If $2\leq \Delta(K_n[M])\leq
\frac{k}{2}$ and $|M|\geq k$, then $\kappa_k(G)\leq
\lambda_k(G)<n-\frac{k}{2}-1$ by Observation \ref{obs1} and Lemma
\ref{lem7}, a contradiction. Therefore, $1\leq |M|\leq k-1$. If
$\Delta(K_n[M])=1$, then $1\leq |M|\leq k-1$ by Lemma \ref{lem8}. We
conclude that $1\leq \Delta(K_n[M])\leq \frac{k}{2}$ and $1\leq
|M|\leq k-1$, as desired.

Conversely, let $G=K_n-M$ where $M\subseteq E(K_n)$ such that $1\leq
\Delta(K_n[M])\leq \frac{k}{2}$ and $1\leq |M|\leq k-1$. In fact, we
only need to show that $\kappa_k(G)\geq n-\frac{k}{2}-1$ for
$\Delta(K_n[M])=1$ and $|M|=k-1$, or $2\leq \Delta(K_n[M])\leq
\frac{k}{2}$ and $|M|=k-1$. The results follow by $(1)$ of Lemma
\ref{lem10} and Lemma \ref{lem11}.\qed\vskip 0.5em

\noindent \textbf{Proof of Theorem \ref{th3}}. If $G$ is a connected
graph satisfying condition $(2)$, then $\kappa_k(G)=n-\frac{k}{2}-1$
by Theorem \ref{th2}. From Observation \ref{obs1}, $\lambda_k(G)\geq
\kappa_k(G)=n-\frac{k}{2}-1$. From this together with Lemma
\ref{lem4}, we have $\lambda_k(G)=n-\frac{k}{2}-1$. Assume that $G$
is a connected graph satisfying condition $(1)$. We only need to
show that $\lambda_k(G)=n-\frac{k}{2}-1$ for $|M|=\lfloor
\frac{n}{2} \rfloor$. The result follows by $(2)$ of Lemma
\ref{lem10} and Lemma \ref{lem4}.

Conversely, assume that $\lambda_k(G)=n-\frac{k}{2}-1$. Since $G$ of
order $n$ is connected, we can consider $G$ as a graph obtained from
a complete graph $K_n$ by deleting some edges. From Corollary
\ref{cor1}, $G=K_n-M$ such that $\Delta(K_n[M])\leq \frac{k}{2}$,
where $M\subseteq E(K_n)$. Combining this with Lemma \ref{lem4}, we
have $|M|\geq 1$ and $\Delta(K_n[M])\geq 1$. So $1\leq
\Delta(K_n[M])\leq \frac{k}{2}$. It is clear that if
$\Delta(K_n[M])=1$ then $1\leq |M|\leq \lfloor\frac{n}{2}\rfloor$.
If $2\leq \Delta(K_n[M])\leq \frac{k}{2}$, then $1\leq |M|\leq k-1$
by Lemma \ref{lem7}. So $(1)$ or $(2)$ holds.\qed\vskip 0.5em

\noindent \textbf{Remark 3.} As we know, $\lambda(G)=n-2$ if and
only if $G=K_n-M$ such that $\Delta(K_n[M])=1$ and $1\leq |M|\leq
\lfloor\frac{n}{2}\rfloor$, where $M\subseteq E(K_n)$. So we can
restate the above conclusion as follows: $\lambda_2(G)=n-2$ if and
only if $G=K_n-M$ such that $\Delta(K_n[M])=1$ and $1\leq |M|\leq
\lfloor\frac{n}{2}\rfloor$, where $M\subseteq E(K_n)$. This means
that $4\leq k\leq n$ in Theorem \ref{th3} can be replaced by $2\leq
k\leq n$.\\

\noindent{\bf Acknowledgement:} The authors are very
grateful to the referees for valuable comments and suggestions, which helped to
improve the presentation of the paper.

\end{document}